\documentclass[10pt]{article}

\usepackage{mathptmx}
\usepackage{amsmath}
\usepackage{amscd}
\usepackage{amssymb}
\usepackage{amsthm}
\usepackage{xspace}
\usepackage[all,tips]{xy}
\usepackage[dvips]{graphicx}
\usepackage{verbatim}
\usepackage{syntonly}
\usepackage{hyperref}
\usepackage[capitalise]{cleveref}
\usepackage{xcolor}
\usepackage{mathpazo}
\usepackage[T1]{fontenc}
\usepackage[utf8]{inputenc} 
\usepackage{lmodern}
\usepackage[numbers,sort&compress]{natbib}

\theoremstyle{plain}
\newtheorem{thm}{Theorem}[section]
\newtheorem{cor}[thm]{Corollary}
\newtheorem{prop}[thm]{Proposition}
\newtheorem{lemma}[thm]{Lemma}
\newtheorem{ques}{Question}

\theoremstyle{definition}
\newtheorem{definition}[thm]{Definition}
\newtheorem*{rem}{Remark}

\newtheorem{remark}{Remark}

\newtheorem{conjecture}{Conjecture}

\DeclareMathOperator{\SL}{SL} \DeclareMathOperator{\PSL}{PSL}

 \DeclareMathOperator{\SO}{SO}
 
\DeclareMathOperator{\PU}{PU} 
\DeclareMathOperator{\SU}{SU} 
\DeclareMathOperator{\Sp}{Sp}

\DeclareMathOperator{\ad}{ad}

\DeclareMathOperator{\hd}{hd}

\newcommand{\eps}{\varepsilon}

\newcommand{\til}{\widetilde}

\newcommand{\ga}{\gamma}
\newcommand{\Ga}{\Gamma}

\newcommand{\La}{\Lambda}

\newcommand{\norm}[1]{\left\vert\left\vert #1\right\vert\right\vert}
\newcommand{\inner}[1]{\left\langle #1 \right\rangle }
\newcommand{\set}[1]{\left\{#1\right\}}

\newcommand{\floor}[1]{\left\lfloor #1 \right\rfloor}

\newcommand{\op}{\operatorname}
\newcommand{\supp}{\op{supp}}
\newcommand{\pa}{\partial}

\newcommand{\B}[1]{\ensuremath{\mathbf{#1}}}

\newcommand{\mc}[1]{\ensuremath{\mathcal{#1}}}

\newcommand{\mf}[1]{\ensuremath{\mathfrak{#1}}}

\newcommand{\N}{\ensuremath{\B{N}}}
\newcommand{\Q}{\ensuremath{\B{Q}}}
\newcommand{\R}{\ensuremath{\B{R}}}

\newcommand\rest[1]{\raisebox{-.5ex}{$|$}_{#1}}

\usepackage{fancyhdr}

\fancypagestyle{plain}

\let\oldmarginpar\marginpar
\renewcommand\marginpar[1]{\-\oldmarginpar[\raggedleft\footnotesize #1]%
{\raggedright\footnotesize #1}}

\begin{document}


\title{Homological dimension of discrete subgroups \\ in higher rank simple Lie groups}
\author{Chris Connell, D. B. McReynolds, Shi Wang}
\maketitle


\begin{abstract}
\noindent We investigate the homological dimensions of discrete subgroups of non-compact simple Lie groups using a flow recently introduced by the authors. Using new estimates on the critical index for a broad class of discrete subgroups, we prove that the homological dimension of non-lattice, discrete, Zariski-dense subgroups is bounded above by $n-\alpha r$ where $n$ is the dimension of the associated symmetric space, $r$ is the real rank, and $\alpha \geq \frac{1}{8}$. Under the assumption that the discrete group has regular limit cone, our bound improves to $\frac{5}{6}n+1$. Additionally, we make some conjectures on the possible homological dimensions for non-lattice discrete subgroups and pose a few questions in the more general semisimple setting.  
\end{abstract}

\section{Introduction}


\subsection{Motivation}

The discrete subgroups of non-compact semisimple Lie groups are central figures in many areas of mathematics. Given a non-compact semisimple Lie group $G$ with associated symmetric space $X = X_G = G/K$ and a discrete, finitely generated, torsion-free, Zariski-dense subgroup $\Gamma$, we obtain a locally symmetric manifold $M=X/\Gamma$ which is a $K(\Gamma,1)$--space. There is a long history of results relating the structure of $G$ with the structure of $\Gamma$. When $\Gamma$ is an irreducible,  cocompact lattice and $G$ has real rank at least two, much of the structure of $\Gamma$ is controlled by $G$. These lattices are arithmetic and super-rigid by seminal work of Margulis, and are conjectured to have the congruence subgroup property by Serre. 

Beyond lattices, a major source of higher rank examples comes from higher Teichm\"uller theory and the theory of Anosov representations, beginning with Hitchin's components and Labourie's formulation of Anosov representations, and developed further through the work of Guichard--Wienhard and Kapovich--Leeb--Porti; see \cite{Hitchin92,Labourie06,GuichardWienhard12,KapovichLeebPorti17}. Another major source of interest are thin subgroups of semisimple Lie groups. A \textbf{thin subgroup of $G$} is an infinite index, Zariski-dense subgroup $\Gamma$ of an arithmetic lattice $\Lambda < G$. These subgroups have connections to quantum computing and also have been studied through the lens of expander families and super strong approximation in linear groups (see \cite{Sarnak2014,FuchsMeiriSarnak2014,Bourgain2013,Breuillard2014,BourgainGamburdSarnak2010,KontorovichEtAl2019,BreuillardGreenGuralnickTao}). Discrete subgroups of Lie groups also arise naturally in monodromy representations arising from differential equations  which has a long and storied history going back to the 1800's (see for instance \cite{BajpaiDonaNitsche2025} and \cite{SinghVenkataramana2014} for recent work with regard to thinness). 

Our interest is in the homology of the discrete subgroups $\Gamma < G$ . One basic invariant is the homological dimension. Given a commutative ring $R$ with identity, we define the \textbf{homological dimension of $\Gamma$ over $R$} to be
\[\op{hd}_R (\Ga):=\inf\{s\in \mathbf{Z}_+\,:\, H_i(\Ga;V)=0 \text{ for any } i>s \text{ and any $R\Ga$--module $V$}\}.\]
It follows that $\mathrm{hd}_\mathbf{Z}(\Gamma) \leq \dim(X)$ with equality if and only if $\Gamma$ is a cocompact lattice. When $\Gamma$ is a non-cocompact lattice, Borel--Serre \cite{BorelSerre} proved that $\mathrm{hd}_\mathbf{Z}(\Gamma) = \dim(X) - r_\Q(\Gamma)$ where $r_\Q(\Gamma)$ is the $\Q$--rank of $\Gamma$. Setting $r=r_G$ to be the real rank of $G$, we have $0 \leq r_\Q(\Gamma) \leq r_G$ for any lattice, and $r_\Q(\Gamma) = 0$ if and only if $\Gamma$ is cocompact. 

\begin{rem}
Non-trivial torsion always implies infinite homological dimension. However, passing to a finite index, torsion-free subgroup (which is ensured by Selberg's Lemma for linear $G$), we can use the virtual homological dimension. For this reason, we will assume that $\Gamma$ is torsion-free.     
\end{rem}

\subsubsection{Simple Lie groups}

We believe that the only finitely generated, discrete, torsion-free, Zariski-dense subgroups of non-compact, simple Lie groups with homological dimension in the range of lattices are lattices. 

\begin{conjecture}\label{Conj:1}
Let $G$ be a non-compact, simple Lie group and $\Gamma$ be a discrete, finitely generated, torsion-free, Zariski-dense subgroup. If $\mathrm{hd}_\mathbf{Z}(\Gamma) \geq \dim(X) - r_G+1$, then $\Gamma$ is a lattice.
\end{conjecture}

\begin{conjecture}\label{Conj:2}
Let $G$ be a non-compact, simple Lie group and let $\Gamma \leq G$ be a torsion-free, cocompact lattice. If $\Delta \leq \Gamma$ is a finitely generated, infinite index subgroup, then $\mathrm{hd}_\mathbf{Z}(\Delta) \leq \dim(X) - r_G$. 
\end{conjecture}

It follows trivially that Conjecture \ref{Conj:1} implies Conjecture \ref{Conj:2}. The dimension $\dim(X) - r_G$ is harder to make conjectures on. First, the $3$--dimensional integral Heisenberg group, which is an infinite index subgroup of $\SL(3,\mathbf{Z})$, has homological dimension 3 where $\dim(X) = 5$ and $r_G=2$ but is not Zariski-dense in $\SL(3,\mathbf{R})$. The lattice $\SL(3,\mathbf{Z})$ has $\mathbf{Q}$--rank $2$ and hence has homological dimension $3$ by Borel--Serre. This shows that there are non-lattice subgroups which are not Zariski-dense and with the same homological dimension as a lattice.  

\begin{ques}\label{Q:1}
Let $G$ be a non-compact, simple Lie group with Property T and let $\Gamma \leq G$ be a discrete, finitely generated, torsion-free, Zariski-dense subgroup. 
\begin{itemize}
\item[(A)]
Are the following equivalent:
\begin{itemize}
\item[(i)]
$\Gamma$ is a lattice.
\item[(ii)]
$\mathrm{hd}_\mathbf{Z}(\Gamma) \geq \dim(X) - r_G$.
\end{itemize}
\item[(B)]
If we further assume that $\Gamma$ contains no non-trivial parabolic elements, does (ii) imply (i)?
\end{itemize}
\end{ques}

\begin{remark}
Among the closed, aspherical, geometric $3$--manifolds, only hyperbolic ones can have a discrete and faithful representation with Zariski-dense image (see also \cite{DoubaTsouvalas2023} for a discussion on graph manifolds).
Whether or not there exists a discrete, faithful representation of a finite volume hyperbolic $3$--manifold group into $\SL(3,\mathbf{R})$ with Zariski-dense image is Question 4.4 in Long--Reid \cite{LongReid} which they attribute to Labourie (see also \cite{LongReid2016} and \cite{BallasLong2020}).  Conjecture \ref{Conj:1} would not exclude this. An affirmative answer to Question \ref{Q:1} (A) would exclude compact ones. This question is also related in spirit to the broader problem of distinguishing flexible higher rank representations from Anosov or relatively Anosov phenomena; see \cite{GuichardWienhard12,KapovichLeebPorti17,GueritaudGuichardKasselWienhard17, CanaryZhangZimmer21}. 
\end{remark}

Up to isogeny, there are two infinite families of non-compact, simple Lie groups without Property T, and they are $\SO(n,1)$ with $n \geq 2$ and $\SU(n,1)$ with $n \geq 2$. For $G =\SO(n,1)$ with $n \geq 2$, Question \ref{Q:1} (A) and (B) have negative answers via bending deformations; see Johnson--Millson \cite{JohnsonMillson}. For $G=\SU(n,1)$, the bending construction has no direct analogue: complex hyperbolic space has no totally geodesic hypersurfaces, and Goldman--Millson \cite{GoldmanMillson} proved that lattices in $\SU(n-1,1)$ are locally rigid in $\SU(n,1)$ for $n \geq 3$.  

\begin{ques}\label{Q:1Complex}
For $n \geq 3$, does there exist a non-lattice, discrete, finitely generated, torsion-free, Zariski-dense subgroup $\Gamma < \SU(n,1)$ with $\mathrm{hd}_{\mathbf{Z}}(\Gamma) = 2n-1$? \end{ques}

\begin{rem}
Question \ref{Q:1Complex} is false for $n=2$. Schwartz \cite{Schwartz} produced a discrete, faithful representation of a closed hyperbolic $3$--manifold group into $\PU(2,1)$ with Zariski-dense image (see also \cite{SchwartzBook}).
\end{rem}
 
Misha Kapovich \cite{Kapovich} proved $\mathrm{hd}_R(\Gamma) \leq \delta(\Gamma) + 1$ when $G = \SO(n,1)$ and $\Gamma$ contains no non-trivial parabolic elements where $\delta(\Gamma)$ is the critical exponent. The first two authors with Farb \cite{CFM16} extended Kapovich's work to the other real rank one Lie groups $\SU(n,1)$ with $n\geq 2$, $\Sp(n,1)$ with $n\geq 2$, and $\mathrm{F}_4^{-20}$. This gave an affirmative answer to Question \ref{Q:1} (B) for the Lie groups $\Sp(n,1)$ with $n \geq 2$ and $\mathrm{F}_4^{-20}$. Specifically, \cite{CFM16} proved that for $G = \Sp(n,1)$ with $n\geq 2$ or $G = \mathrm{F}_4^{-20}$, if $\Gamma < G$  is a discrete, torsion-free subgroup with no non-trivial parabolic elements, then 
\[ \mathrm{hd}_R(\Gamma) \leq \dim(X) -2 \text{ or }  \mathrm{hd}_R(\Gamma) \leq \dim(X) - 4, \] 
respectively.  These are precisely the non-compact, simple Lie groups with real rank one and Property T. The homological gaps are directly connected to explicit Kazhdan constants via  Corlette \cite{Corlette} who proved a gap theorem for critical exponents of discrete, non-lattice subgroups of these Lie groups. The main upgrade in \cite{CMW23} was replacing the natural map employed in \cite{Kapovich} and \cite{CFM16} with the natural flow. However, all of these methods use Gromov's Mass Vanishing theorem (\cite{Kapovich} or \cite[App]{CMW23}) and the $k$--volume contracting behavior of either the natural map or the natural flow to establish homological vanishing. This was the key insight in \cite{Kapovich}.

\subsubsection{Semisimple Lie groups}

For semisimple Lie groups, we need to take additional care in how we state our questions. First, one should assume that $G$ is isotypic to ensure that it admits irreducible, cocompact lattices. Given the role Property T played in \cite{CFM16}, we will also assume that each simple factor has rank at least two to ensure that the semisimple Lie group has Property T. We also wish to exclude (virtual) products as examples and so will impose conditions on the projection maps for this purpose. 

\begin{ques}\label{Q:2}
Let $G$ be an isotypic, semisimple Lie group with non-compact simple factors $G_1,\dots,G_k$ and $r_{G_j}>1$ for $j=1,\dots,k$. If $\Gamma < G$ is a discrete, finitely generated, torsion-free, Zariski-dense subgroup which projects densely onto each simple factor $G_j$ of $G$ and $\mathrm{hd}_\mathbf{Z}(\Gamma) \geq \dim(X) - r_G$, then is $\Gamma$ an irreducible lattice?
\end{ques}

We believe the gaps in homological dimension should be controlled by the real rank and not by Property T. The simplest example is $G = (\PSL(2,\mathbf{R}))^2$ where $\dim(X) = 4$ and $r_G = 2$. Closed $\mathrm{Sol}$ 3--manifold groups have discrete, faithful representations to $(\PSL(2,\mathbf{R}))^2$ as cusp groups of Hilbert modular surface groups which are irreducible, non-cocompact lattices in $G$. This provides examples of discrete subgroups with homological dimension $3$ which are not lattices; Hilbert modular surface groups have $\mathbf{Q}$--rank 1 and so also have homological dimension 3 by Borel--Serre. However, as before, these cusp subgroups are not Zariski-dense. Additionally, given a pair of distinct hyperbolic structures on a genus $g>1$ closed, orientable surface $\Sigma_g$, we have a pair of non-conjugate, discrete, faithful representations $\rho_1,\rho_2\colon \pi_1(\Sigma_g) \to \PSL(2,\mathbf{R}).$ The image of $\rho_1 \times \rho_2$ is a discrete, Zariski-dense subgroup of $\PSL(2,\mathbf{R}) \times \PSL(2,\mathbf{R})$ and has homological dimension 2. Note in this case, the image projects discretely to each $\PSL(2,\mathbf{R})$ factor. Note additionally that Teichm\"{u}ller curves also can produce discrete examples of Fuchsian groups (see \cite{McMullen} and \cite{McMullen25}).    

\begin{ques}\label{Q:3}
Let $G$ be an isotypic, semisimple Lie group with non-compact simple factors $G_1,\dots,G_k$ and $k>1$. If $\Gamma < G$ is a discrete, finitely generated, torsion-free, Zariski-dense subgroup which projects densely onto each simple factor $G_j$ of $G$ and $\mathrm{hd}_\mathbf{Z}(\Gamma) \geq \dim(X) - r_G$, then is $\Gamma$ an irreducible lattice?
\end{ques}

Natural classes of examples to explore for Question \ref{Q:3} are images of representations of non super-rigid lattices into semisimple Lie groups, of which free and surface groups are the most basic examples. In connection with the limit cone theory of Benoist and Quint \cite{Benoist97,Quint-divergence02,Quint02,Quint03}, Benoist asked if there exists a discrete, finitely generated free subgroup in $(\SL(2,\mathbf{R}))^2$ which is Zariski-dense and has dense projection onto each $\SL(2,\mathbf{R})$ factor (see \cite[Remark below Cor 3.3]{BrodyEtAl25} and the related questions in  \cite[Questions 1 and 2]{FisherEtAl18}). This would not be excluded if Question \ref{Q:3} had an affirmative answer. However, it would exclude discrete, Zariski-dense surface subgroups of $(\SL(2,\mathbf{R}))^2$ which have dense projections onto each factor; note that having dense projections on each factor imposes restrictions on the Toledo invariant of the representation (see \cite{BurgerIozziWienhard}). For higher dimensional real hyperbolic examples, an affirmative answer to Question \ref{Q:3} would not allow for discrete, Zariski-dense closed hyperbolic $3$--manifold groups in $\PSL(2,\mathbf{R}) \times \PSL(2,\mathbf{C})$ which have dense projections onto each factor but would for $(\PSL(2,\mathbf{R}))^2 \times \PSL(2,\mathbf{C})$. One could also consider other rank one lattices that are not super-rigid. However, by work of Carlson--Toledo \cite{CarlsonToledo}, if $\Gamma$ is a cocompact lattice in $\SU(n,1)$ with $n>1$, then any discrete representation of $\Gamma \to \PSL(2,\mathbf{R})$ either has cyclic image or factors through a homomorphism $\Gamma \to \Delta \to \PSL(2,\mathbf{R})$ where $\Delta$ is a surface group. It follows that there can exist no discrete and faithful representation of $\Gamma \to (\PSL(2,\mathbf{R}))^k$ for any $k>0$. The isotypic semisimple Lie groups of type $\mathrm{A}_1$, like the ones above, are a good place to start looking for examples that might show Question \ref{Q:3} has a negative answer. For higher dimensional analogs, one might hope to find a closed hyperbolic $n$--manifold group $\Gamma$ with a rich representation theory to $ \SO(n,1) \times \SO(n,1)$ to find a discrete, faithful representation with Zariski-dense image and dense projections onto each factor. However, in this case $\dim(X) = 2n$, the image would have $\mathrm{hd}(\Gamma) = n$, and $r=2$. For complex hyperbolic, we would instead want to find a lattice in $\SU(n,1)$ with a discrete, faithful and Zariski dense representation to  $\SU(n,1) \times \SU(n,1)$. Here $\dim(X) = 4n$, $\mathrm{hd}(\Gamma) = 2n$, and $r=2$. In particular, we only get $\mathrm{hd}(\Gamma) = \dim(X) - r$ when $G=\SO(2,1)$, $\SU(1,1)$, and $\PSL(2,\mathbf{R})$ which are all isogenous.         

\subsection{Main results}

Our main result establishes a weaker version of Conjecture \ref{Conj:1} under a mild condition on $\Gamma$ which we now explain. In our recent paper \cite{CMW23}, we introduced a flow on $M$ inspired by the natural map of Besson--Courtois--Gallot \cite{BCG} where we proved that in dimensions $k \geq j_X(\Gamma)$, the $k$--volume of the flow is contracting. Using this flow, we proved the following homological vanishing result:

\begin{thm}\cite{CMW23}\label{thm:Hadamard}
If $X$ is a Hadamard space with bounded sectional curvatures, $\Gamma < \mathrm{Isom}(X)$ is a discrete, torsion-free subgroup with $M = X/\Gamma$, and $\mathbf{V}$ is a flat bundle over $M$, then for every $\epsilon > 0$, the homomorphism $i_k:H_k(M_{\eps};\mathbf{V})\rightarrow H_k(M;\mathbf{V})$ induced by inclusion is surjective whenever $k\geq  j_X(\Gamma)$.
\end{thm}

Here $M_\eps$ is the $\eps$--thin part of $M$, and $j_X(\Gamma)$ is the critical index of $\Gamma$, given by the following.

\begin{definition}\label{D:CritIndex} 
Let $X$ be an $n$--dimensional Hadamard manifold, $\Ga<\op{Isom}(X)$ a discrete, torsion-free subgroup, and $\mu_x$ be any conformal density of dimension $\delta_\mu$. We define the \textbf{critical index of $\Ga$ associated to $\mu$} to be
\[j_X(\Gamma,\mu) :=\min\left(\{k\in \mathbf{N}: \inf_{(x,\theta)\in X\times \partial_\mu X} \op{tr}_k(\nabla dB_{(x,\theta)})>\delta_\mu\}\cup\set{n+1}\right).\]	
and the \textbf{critical index} of $\Ga$
\[ j_X(\Gamma)=\inf \set{j_X(\Gamma,\mu): \mu \text{ is a } \delta_\mu-\text{conformal density}}.\]
\end{definition}

Here $\nabla dB_{(x,\theta)}$ denotes the Hessian of the Busemann function at $x$ towards $\theta$, $\partial_\mu X$ is the support of $\mu$, and $\op{tr}_k$ is the sum of first $k$ smallest eigenvalues. We recall the definition of a conformal density of dimension $\delta$ for $\Gamma$ (also known more simply as a $\delta$--conformal density).
\begin{definition}\label{def:delta_conf_density}
Given $\delta> 0$, we say that a family of finite, nonzero, Borel measures $\{\mu_x:x\in X\}$ is a \textbf{$\delta$--conformal density on $X$} (associated to the $\Gamma$--action) if
\begin{enumerate}
\item each $\mu_x$ is supported on $\supp{\mu_x}\subset\partial_\infty X$,
\item $\mu_x$ is $\Ga$--equivariant, that is, $\mu_{\ga x}=\ga_*\mu_x$ for all $x\in X$ and $\ga\in \Ga$,
\item the measures $\mu_x$ are pairwise equivalent with Radon--Nikodym derivatives given by
\[\frac{d\mu_x}{d\mu_y}(\theta)=e^{-\delta B_\theta(x,y)}\]
for all $x,y\in X$ and $\theta\in \supp{\mu_x}$.
\end{enumerate}
\end{definition}

\subsubsection{Critical index bounds}

The most difficult aspect of employing this flow to obtain homological vanishing results is estimating $j_X(\Gamma)$. In this paper, we provide new estimates for $j_X(\Gamma)$ under certain conditions on $\Gamma$ arising in higher rank Patterson--Sullivan theory. As the contracting nature of the flow is of broader use, we first state our results on $j_X(\Gamma)$. We restrict our context and assume $X=G/K$ is an irreducible, higher rank symmetric space of non-compact type and $\Gamma<G$ is a discrete, torsion-free subgroup. Using the general scheme of \cite[Thm 1.11]{CMW23} and \cite[Prop 5.27]{Fraczyk-Lowe}, we give a sharpened estimate on $j_X(\Gamma)$ for all irreducible higher rank symmetric spaces. The condition we require on $\Gamma$ is the following (we refer the reader to \S 2 for the notation):
\begin{equation}\label{eq:cond}
\textrm{ There exists }\lambda\in (0,1] \textrm{ such that }\lambda(2\rho-\Theta)\textrm{ extremizes }\psi_\Gamma\textrm{ at a point }z \in\mathfrak a^{++}.
\end{equation}

\begin{remark}
Condition (\ref{eq:cond}), roughly speaking, means that the $\Ga$--orbits, when measured via the $(2\rho-\Theta)$ distance, dominate the growth in a regular direction.
\end{remark}

\begin{thm}\label{thm:j_X_bound}
Let $G$ be a non-compact, real, simple Lie group with real rank $r_G>1$, $n_G=\dim(G/K)$, and $\Gamma<G$ a discrete, finitely generated, torsion-free, Zariski-dense subgroup which is not a lattice. Suppose further Condition \eqref{eq:cond} holds. Then we have the following upper bound on $j_X(\Gamma)$ hold, according to the type of restricted root system of $G$, where $r=r_G$.  
\begin{enumerate}
\item If $G$ is of Type $\mathrm{A}_r$, then $j_X(\Gamma)< n_G+1-\frac{r}{8}$.
\item If $G$ is of Type $\mathrm{B}_r$ or $\mathrm{D}_r$, then $j_X(\Gamma)< n_G+1-\frac{3r}{16}$.
\item If $G$ is of Type $\mathrm{C}_r$ or $(\mathrm{BC})_r$, then $j_X(\Gamma)< n_G+1-\frac{r}{4}$.
\end{enumerate}
\end{thm}

These results show that the natural flow is volume contracting in codimensions proportionate to $r_G$. Under a more restrictive hypothesis on the limit cone $\mathcal{C}_\Gamma$ of $\Gamma$ (see \S 2 for the definition of $\mathcal{C}_\Gamma$), we show the flow is contracting in codimensions proportionate to $n_G$. One class of examples that satisfies this condition on the limit cone are Borel--Anosov subgroups. Anosov representations were introduced by Labourie for surface groups \cite{Labourie06}, Guichard--Wienhard extended the theory to word-hyperbolic groups and developed domains of discontinuity \cite{GuichardWienhard12}, and Kapovich--Leeb--Porti gave dynamical and coarse-geometric characterizations \cite{KapovichLeebPorti17}. 

\begin{thm}\label{thm:anosov}
Let $G$ be a non-compact, real, simple Lie group with $r_G>1$ and which is not of type $\mathrm{A}_2$ restricted root system, $n_G = \dim(G/K)$, and $\Gamma<G$ be a discrete, finitely generated, torsion-free, Zariski-dense subgroup with regular limit cone (i.e. $\mathcal C_\Gamma\subset \mathfrak a^{++}\cup\set{0}$). Then we have  
\[ j_X(\Gamma)< \frac{5n_G}{6}+1. \]
In particular, if $\Gamma$ is a Borel--Anosov subgroup of $G$, then this holds.
\end{thm}

There is by now a substantial parallel literature on the geometry and dynamics of Anosov subgroups, including proper actions and tameness of locally symmetric quotients \cite{GueritaudGuichardKasselWienhard17,GuichardKasselWienhard18}, pressure metrics and entropy/intersection theory \cite{BridgemanCanaryLabourieSambarino15}, and conformal measures, local mixing, and counting in higher rank \cite{EdwardsLeeOh20,EdwardsLeeOh21,LeeOh22}. Recent work of \cite{Fraczyk-Lowe} employed the natural flow to establish higher dimensional analogs of Property (T). The above estimates provide additional applications to other semisimple Lie groups and stronger fixed point results for groups $\Gamma$ with regular limit cones.

\subsubsection{Homological dimension bounds}

Using these estimates on $j_X(\Gamma)$ and following the methods of \cite{CMW23}, we obtain the following upper bounds on the homological dimension in the above two settings.

\begin{cor}\label{cor:hd_bound}
Let $X, \Gamma, n_G, r_G$ be as in Theorem \ref{thm:j_X_bound}. Suppose the injectivity radius of $\Gamma\backslash X$ is uniformly bounded from zero. Below, $r=r_G$.
\begin{enumerate}
\item If $G=\SL(r+1,\mathbf{R}), \SL(r+1,\mathbf{C})\textrm{ or } \SU^*(2r+2)$, then $\hd_R(\Gamma)< n_G-\frac{r}{8}$.
\item If $G=\SO(r,s)$ with $r\leq s$, $\SO(2r,\mathbf{C})$ or $\SO(2r+1,\mathbf{C})$, then $\hd_R(\Gamma)< n_G-\frac{3r}{16}$.
\item If $G=\Sp(r,\mathbf{R})$, $\Sp(r,\mathbf{C})$, $\SU(r,s)$ with $r\leq s$, $\Sp(r,s)$ with $r\leq s$, $\SO^*(4r)$ or $\SO^*(4r+2)$, then $\hd_R(\Gamma)< n_G -\frac{r}{4}$.
\end{enumerate}
\end{cor}

This vanishing result is similar in spirit to Matsushima's vanishing theorem \cite{Matsushima} for the rational homology of cocompact lattices $\Gamma$ in semisimple Lie groups $G$. Specifically, Matsushima proved $H_{n-k}(\Gamma,\Q) = 0$ for $0 < k < r$ where $n=\dim(G/K)$ and $r$ is the real rank of $G$.

\begin{remark}
We believe that the condition on the injectivity radius can be removed. This amounts to saying that the homological dimension contributed by the ``zero thin at infinity'' is bounded by $n_G-r_G$, which requires further investigation on the structure of the thin part for general infinite volume locally symmetric spaces.
\end{remark}

We also obtain stronger bounds on the homological dimension with our hypothesis on the limit cone of $\Gamma$.   

\begin{cor}
Let $G$ be a non-compact, real, simple Lie group with $r_G>1$ and which is not of type $\mathrm{A}_2$ restricted root system, $n_G = \dim(G/K)$, and $\Gamma<G$ be a discrete, finitely generated, torsion-free, Zariski-dense subgroup with regular limit cone (i.e. $\mathcal C_\Ga\subset \mathfrak a^{++}\cup\set{0}$). Suppose the injectivity radius of $\Gamma\backslash X$ is uniformly bounded from zero. Then we have  
\[ \mathrm{hd}_R(\Gamma)< \frac{5n_G}{6}+1. \]
In particular, if $\Ga$ is a Borel--Anosov subgroup of $G$, then this holds.
\end{cor} 

In the Anosov case of the preceding corollary, Canary--Tsouvalas in  \cite{CanaryTsouvalas20} proved that if a torsion-free hyperbolic group admits a $P_k$--Anosov representation into $\SL(d,\mathbf{R})$, then $\operatorname{cd}(\Gamma)\le d-k$, hence $\operatorname{hd}_{\mathbf{Z}}(\Gamma)\le d-k$, except in the four Hopf-fibration cases $(d,k)=(2,1),(4,2),(8,4),(16,8)$, where the bound is $\operatorname{cd}(\Gamma)\le d-k+1$. These bounds are better than the ones given above. However, the condition that the limit cone be regular is less restrictive than being a Borel--Anosov subgroup. One cannot expect the same vanishing results for this larger class, though these vanishing results are substantially better than the much more general condition (\ref{eq:cond}). 

The Anosov case should be compared with several recent lines of work on quantitative and dynamical features of higher rank representations. Bridgeman--Canary--Labourie--Sambarino developed the pressure metric along with the entropy/intersection theory for projective Anosov representations \cite{BridgemanCanaryLabourieSambarino15}; Pozzetti--Sambarino--Wienhard obtained rigidity results for Anosov representations with Lipschitz limit set \cite{PozzettiSambarinoWienhard23}; and Canary--Zhang--Zimmer extended parts of the Anosov/Hitchin picture to geometrically finite Fuchsian groups \cite{CanaryZhangZimmer21,CanaryZhangZimmer22}. These works study growth, entropy, and limit set regularity, while we use higher rank Patterson--Sullivan theory and the natural flow to obtain homological dimension bounds.

\subsection{Some final questions}

We end the introduction with questions related more directly to this paper. For Anosov representations, several works show that dynamical data such as periods, entropy, pressure, and growth indicators vary in highly structured ways \cite{Sambarino14,Sambarino14Counting,BridgemanCanaryLabourieSambarino15,KimEtAl21,PozzettiSambarinoWienhard23}. Question \ref{Q:4} asks whether the topological invariant $\mathrm{hd}_R(\Gamma)$ can be made comparably efficient with respect to the critical index threshold. Assuming $\Gamma < G$ is a discrete, finitely generated, torsion-free, Zariski-dense subgroup, we know that $j_X(\Gamma) - \mathrm{hd}_R(\Gamma) > 0$. It could be that $j_X(\Gamma)$ is much larger than $\mathrm{hd}_R(\Gamma)$ (e.g. a discrete free group in $\SO(n,1)$ with larger critical exponent). We ask if there exists an efficient representation with regard to witnessing the homological dimension of $\Gamma$.

\begin{ques}\label{Q:4}
Let $G$ be a non-compact semisimple Lie group without compact factors and let $\Gamma<G$ be a discrete, finitely generated, torsion-free, Zariski-dense subgroup. Does there exist a non-compact semisimple subgroup $H \leq G$ and faithful representation $\rho\colon \Gamma \to H$ with discrete, Zariski-dense image such that $0 < j_{X_H}(\rho(\Gamma)) - \mathrm{hd}_R(\Gamma) \leq 1?$
\end{ques}

We do not know of a single example of a Zariski-dense subgroup of a non-compact simple Lie group which does not satisfy condition (\ref{eq:cond}).

\begin{ques}\label{Q:5}
Let $G$ be a non-compact simple Lie group.
\begin{itemize}
\item[(A)] If $\Gamma$ is a cocompact lattice in $G$ and $\Delta < \Gamma$ is an infinite index subgroup, does $\Delta$ satisfy (\ref{eq:cond})?
\item[(B)]
Does there exist a discrete, finitely generated, Zariski-dense subgroup $\Gamma < G$ which does not satisfy (\ref{eq:cond})?\end{itemize}
\end{ques}

The homological dimension is non-increasing in passing to subgroups. We believe the same should hold for the critical index.

\begin{ques}\label{Q:6}
Let $G$ be a non-compact simple Lie group. If $\Gamma < G$ is a discrete, finitely generated, torsion-free, Zariski-dense subgroup and $\Delta < \Gamma$ is a finitely generated subgroup, do we have $j_X(\Delta) \leq j_X(\Gamma)$?
\end{ques}

Our final question is the analog of Conjecture \ref{Conj:1} but for the critical index.

\begin{ques}\label{C:7}
Let $G$ be a non-compact simple Lie group and let $\Gamma$ be a discrete, finitely generated, torsion-free, Zariski-dense subgroup. If $\Gamma$ is not a lattice, is $j_X(\Gamma) < \dim(X) - r_G + 1$?
\end{ques}

\subsection*{Acknowledgments.}

This work was partially supported by grants Simons Foundation/SFARI-965245 (CC) and NSF DMS-1812153 (DBM) and NSFC No.12301085 (SW). We would like to thank Subhadip Dey, Alan Reid, Matthew Stover, and Hongwei Zhang for helpful discussions. We would especially like to thank Ben Lowe for pointing out the possible improvement of our previous work.

\subsection*{Misha Kapovich}

Near the completion of this work, we learned of the passing of Misha Kapovich. This paper does not exist without the elegant ideas of Misha from \cite{Kapovich}. We will sadly miss him as a friend. We are thankful that he will live on in the work of generations of mathematicians who were truly inspired by him, like us, to go further into the unknown. 

\section{Patterson--Sullivan theory in higher rank}\label{sec:application}

In this section, we apply Theorem \ref{thm:Hadamard} from \cite{CMW23} to the setting of irreducible, higher rank, locally symmetric spaces. We note that, as opposed to the case of rank one symmetric spaces, the estimate of the critical index $j_X(\Gamma)$ becomes significantly more difficult, partly due to the following two reasons, which make the choice of $\theta$ strategic.

\begin{itemize}
\item[$\bullet$] The Hessian of the Busemann function $\nabla dB_{(x,\theta)}$ now depends on the $\theta$ direction.
\item[$\bullet$] The existence of the $\delta$--conformal density follows from the deep theory (the extension of Patterson--Sullivan theory) formulated by Benoist and Quint \cite{Benoist97,Quint-divergence02,Quint02,Quint03} (see also \cite{Link04}). For more recent developments on conformal densities, local mixing, counting, and temperedness in higher rank, see \cite{EdwardsOh22,EdwardsLeeOh20,EdwardsLeeOh21,LeeOh22}. We apply it with a special choice of the $G$--orbit of $\theta$, and in particular $\delta$ will depend on such choice.
\end{itemize}

\subsection{Setup}

Let $G$ be a semisimple Lie group of noncompact type which we will assume has real rank at least two, and let $\Ga<G$ be a discrete, torsion-free, Zariski-dense subgroup. Denote the associated symmetric space by $X=G/K$ where $K<G$ is a maximal compact subgroup. We denote the quotient locally symmetric manifold by $M=\Ga\backslash X$. For the following, we choose a basepoint ${p}\in X=\til{M}$ once and for all. We follow the conventional notation of letting $A<G$ denote a maximal diagonalizable, abelian subgroup with Lie algebra $\mf{a}$. Denote by $\Sigma\subset \mf{a}^*$ the set of all roots in the restricted root system, and $\Sigma^+\subset\Sigma$ a set of positive roots with respect to a choice of positive system, and $\Pi\subset \Sigma^+$ the corresponding set of simple roots. This determines a positive \textbf{Weyl chamber} 
\[ \mf{a}^+=\bigcap_{\alpha\in \Pi} \set{H\in \mf{a}\,:\, \alpha(H)\geq 0}. \]
For $\theta\subset \Pi$, we set
\[ \mf{a}_\theta=\bigcap_{\alpha\in \Pi \setminus \theta} \ker \alpha, \quad \mf{a}_\theta^+=\mf{a}^+\cap \mf{a}_\theta\quad\text{and}\quad \mf{a}_\theta^{++}=\mf{a}_\theta^+\setminus\left(\bigcup_{\tau\subsetneq \theta} \mf{a}_{\tau}^+\right).\]

Note that the subsets of $\Pi$ are in bijective correspondence with the faces of the canonical Weyl chamber $\mf{a}^+$ associated to $\Sigma^+$. When $\theta=\Pi$ we drop the subscript writing $\mf{a}^{++}$ for $\mf{a}^{++}_\Pi$.

\subsection{The Benoist Cone}\label{sec:Benoist_Cone}

We fix a Weyl chamber $\mf{a}^+$ in $\mf{a}$. For any $g\in G$, let $a(g)$ be the unique element in $\mf{a}^+$ such that 
\[ g=k_1(g) (\exp a(g)) k_2(g) \] 
for $k_1(g),k_2(g)\in K$ that is guaranteed by the Cartan decomposition $G=KA^+K$. (The element $k_1(g)$ is only uniquely defined if $a(g)\in \mf{a}^{++}$.) The map $a\colon G\to \mf{a}^+$ is sometimes called the \textbf{Cartan projection}.  Define the \textbf{limit cone} $\mc{C}_\Ga$ of $\Ga$ to be the set of $x\in \mf{a}^+$ such that there exist sequences $t_n>0$ with $t_n\to 0$ and $\ga_n\in\Ga$ such that $t_n a(\ga_n)\to x$. Benoist \cite{Benoist97} proved that this cone is convex and has nonempty interior. Related constructions of discrete, Zariski-dense subgroups and convex projective examples also appear in \cite{BenoistConvexesI,BenoistConvexesII,BenoistConvexesIII,BenoistConvexesIV}. It is immediate from the definition that $\La(\Gamma)\subset K\cdot \pa_\infty \mc{C}_\Ga$ where $\Lambda(\Gamma)$ is the limit set of $\Gamma$. However, the reverse inclusion may not hold.

We let $\theta_\Gamma\subset \Pi$ denote the subset corresponding to the smallest closed face of the Weyl chamber at infinity whose $G$--orbit contains the limit set $\Lambda(\Gamma)$.  In \S 7 of \cite{Benoist97} it was shown that $\theta_\Ga=\Pi$ for any reductive group defined over $\R$.  Fix a norm $\norm{\cdot}$ on $\mf{a}$.  For every open cone $\mc{C}$ based at $0$ in $\mf{a}$, let $\delta_{\mc{C}}(\Gamma)\in \set{-\infty}\cup[0,\infty)$ be the \textbf{critical exponent} of convergence of the series
\[ \sum_{\substack{\ga\in\Gamma \\ a(\ga)\in \mc{C}}} e^{-s \norm{a(\ga)} }.\]
(Observe that if $\mc{C}$ does not intersect $\mf{a}^+$, then $\delta_{\mc{C}}=-\infty$ since $a(\ga)\in \mf{a}^+$ for all $\ga\in \Gamma$.) We also write $\delta_{\norm{\cdot}}(\Gamma)$ for $\delta_{\mf{a}}(\Gamma)=\delta_{\mf{a}^+}(\Gamma)$, that is the critical exponent for $\Gamma$ with respect to $\norm{\cdot}$ for the full cone $\mc{C}=\mf{a}^+$. Following Benoist and Quint \cite{Benoist97,Quint-divergence02,Quint02,Quint03}, we define a function $\psi_\Gamma\colon \mf{a}\to \mathbf{R} \cup \set{-\infty}$, called the \textbf{ $\Gamma$--indicator function}, by setting $\psi_\Ga(0)=0$ and for $x\neq 0$ we set
\begin{align}\label{eq:Gamma-indicator}
\psi_\Ga(x)=\norm{x} \inf_{\mc{C}}\delta_{\mc{C}}(\Ga),
\end{align}
where the infimum is over all open cones $\mathcal{C}$ based at $0\in \mf{a}^+$ containing $x$. In \cite{Quint-divergence02} it was shown that $\psi_\Ga$ does not depend on the choice of norm $\norm{\cdot}$ and that it is concave and upper semicontinuous. Moreover, Quint shows that $\psi_\Ga(x)\geq 0$ on $\mc{C}_\Gamma$, $\psi_\Ga(x)> 0$ on the interior of $\mc{C}_\Gamma$, and $\psi_\Ga(x)=-\infty$ everywhere else. For the purpose of this paper, we choose the norm that comes from the symmetric Riemannian inner product. 

\subsection{Generalized Busemann Functions}

Given  $g\in G$ we write for $g=k_ga_gn_g$ for the corresponding projections $k_g\in K$, $a_g\in A$ and $n_g\in N$ of the element $g$ in the  Iwasawa decomposition $KAN$ of $G$. We define $f\colon G \to \mathfrak{a}$ to be given by $f(g) \in \mathfrak{a}$ where $\exp(f(g)) = a_g$ or equivalently $f(g) = \log(a_g)$, where $\log\colon A \to \mathfrak{a}$. Observe that $f(g)\in \mf{a}$ is not in general the same as the Cartan projection $a(g)\in \mf{a}^+$.

We define the \textbf{Furstenberg boundary} to be $\mc{F}=G/P$ where $P$ is the minimal parabolic subgroup stabilizing the maximal Weyl chamber at infinity. Observe that by the associated Langlands decomposition, $P=M A N$ where $M<K$ is the stabilizer in $K$ of the face $\exp(\mf{a}^+)\cdot p$ (and hence of all $\exp(\mf{a})\cdot p$).  Hence we may identify the Furstenberg boundary as $\mc{F}=K/M$. 

Define the \textbf{$\mf{a}$--valued Busemann cocycle} $\beta\colon \mc{F} \times G \times G \rightarrow \mathfrak{a}$ as follows. For $[k] \in \mc{F}$ and $g, h \in G$,
\[ \beta_{[k]}(g, h):=f(g^{-1}k)-f(h^{-1}k).\]
Observe that the above definition does not depend on the choice of $k$, and moreover, $\beta_\xi(g,h)$ is right $K$--invariant in $g$ and $h$ separately and is left $G$--equivariant since for $\sigma\in G$ we have 
\begin{align*}
\beta_{\sigma[k]}(\sigma g, \sigma h)&=f(g^{-1}k(\sigma k)^{-1}k_{\sigma k})-f(h^{-1}k(\sigma k)^{-1}k_{\sigma k})\\
&=f(g^{-1}k n_{\sigma k}^{-1}a_{\sigma k}^{-1}k_{\sigma k}^{-1}k_{\sigma k})-f(h^{-1}k n_{\sigma k}^{-1}a_{\sigma k}^{-1}k_{\sigma k}^{-1}k_{\sigma k})\\
&=f(g^{-1}k a_{\sigma k}^{-1} n')-f(h^{-1}k a_{\sigma k}^{-1} n')\\
&=f(g^{-1}k) - f(a_{\sigma k}) -f(h^{-1}k)+f(a_{\sigma k})=\beta_{[k]}( g, h).
\end{align*}
Here we have used that $f$ is right $N$--invariant and that $n'$ is some element of $N$ resulting from the fact that $A$ normalizes $N$. This also implies that $f(g a)=f(g)+f(a)$ for all $g\in G$ and $a\in A$. The $K$--invariance of $\beta_\xi$ for $\xi\in\mc{F}$, implies that we can treat $\beta$ as a function $\beta\colon \mc{F}\times X\times X\to \mf{a}$. We will distinguish both versions of $\beta$ by context. 

We will call the $G$--orbit of the interior of the canonical Weyl chamber at infinity, $\pa_\infty (\exp(\mf{a}^{++})\cdot p)$, the \textbf{regular set at infinity}, and we denote it by $\pa_R X$. Since $\theta_\Gamma = \Pi$, $\Lambda(\Gamma)$ always intersects the regular set at infinity. For any $\xi\in \pa_RX$, we may write $\xi=k\xi_0$ with $\xi_0=\lim_{t\to \infty}\exp(t a_{0})$ for some unit vector $a_{0}\in \mf{a}^{++}$. In particular, the $G$--orbits in $\pa_R X$ coincide with the $K$--orbits which are each homeomorphic to $\mc{F}=K/M=G/P$. For this $\xi$, the ordinary Busemann function relates to the $\mf{a}$--valued Busemann cocycle by $\inner{\beta_{[k]}(x,y),a_{0}}=B_{\xi}(x,y)$ (cf. \cite[Lemma 8.9]{Quint02}).

\subsection{\texorpdfstring{$(\Ga,\phi)$}{}--conformal densities}

Now we consider a generalization of the Patterson--Sullivan construction. First we state some definitions; compare the following definition with the equivalent \cite[Defn 2.2]{EdwardsOh22} and the definition before \cite[Thm 7.2]{Quint02}. (We remark that there is another complete and parallel development of Patterson--Sullivan densities in higher rank spaces due to Link \cite{Link04}, though we will not use that formulation here.) In what follows recall the notations from \cref{sec:Benoist_Cone}.

\begin{definition}
We say that $g_i \to \infty$ \textbf{regularly} if $\alpha(a(g_i))\to \infty$ for all $\alpha\in \Pi$. The \textbf{$\mc{F}$--limit set of $\Gamma$} is the subset of $\mc{F}$ given by,
\[ \Lambda_{\mc{F}}(\Gamma)=\set{[k]\in \mc{F} \,:\, \text{ there exists }\, \ga_n\to\infty \text{ regularly with } \lim_n [k_1(\ga_n)]\to [k]}.\]
\end{definition} 

\begin{remark}\label{rem:limit_sets_same}
By \cite[Thm 3.6]{Benoist97} and \cite[Thm 2.5]{Guivarch90}, $\Lambda_{\mc{F}}(\Gamma)$ is Zariski dense in $\mc{F}$ and the unique $\Ga$--invariant closed subset of $\mc{F}$. In particular it is the closure of any of its $\Gamma$--orbits. As a consequence, the relation to the ordinary limit set $\Lambda(\Gamma)\subset \pa_\infty X$ is rather simple. While from the definition we a priori have,
\[\Lambda_{\mc{F}}(\Gamma)=\bigcup_{\xi\in\pa_\infty \exp(\mf{a}^{++})\cdot p} \set{[k]\in \mc{F}\,:\, k\xi\in \Lambda(\Gamma)},\]
in fact the $G$--strata of $\Lambda(\Gamma)$ are either identified with $\Lambda_{\mc{F}}(\Gamma)$ or empty. Specifically, for each $\xi\in \pa_\infty \exp(\mf{a}^{++})\cdot p$, either $\Lambda(\Gamma)\cap G\cdot\xi=\emptyset$ or else 
\[ \Lambda(\Gamma)\cap G\cdot\xi =\set{k \xi\,:\, [k]\in\Lambda_{\mc{F}}(\Gamma)}. \] 
This follows immediately from the fact that $\set{[k]\in \mc{F}\,:\,k\xi\in\Lambda(\Gamma)}$ is a closed $\Gamma$--minimal subset of $\mc{F}$. (The analogous statement is also true in each singular strata where $\mc{F}_\theta=K/M_\theta$ for $\theta\subset \Pi$ replaces the role of $\mc{F}=\mc{F}_\Pi$, but we will not need this.)
\end{remark}

\begin{definition}
Given $\phi \in \mf{a}^{*}$ and a discrete, Zariski-dense subgroup $\Ga<G$, a Borel probability measure $\nu$ on $\mc{F}$ is called a \textbf{$(\Gamma,\phi)$--conformal measure} if for any $\gamma \in \Gamma$ and $[k] \in \mc{F}$,
\[ \frac{d \gamma_{*} \nu}{d \nu}([k])=e^{-\phi\left(\beta_{[k]}(\gamma,e)\right)}.\]
\end{definition}

Similarly, we can define $(\Ga,\phi)$--conformal densities as follows.

\begin{definition}
A \textbf{$(\Ga,\phi)$--conformal density} is a family of finite nonzero Borel measures $\set{\nu_x}_{x\in X}$ on $\mc{F}$ such that
\begin{itemize}
\item $\frac{d\nu_x}{d\nu_y}([k])=e^{-\phi(\beta_{[k]}(x,y))}$ for all $x,y\in X$ and $[k]\in \mc{F}$
\item $\ga_*\nu_x=\nu_{\ga x}$ for all $\ga\in \Ga$ and $x\in X$.
\end{itemize}
\end{definition}

Notice that if $\set{\nu_x}$ is a $(\Gamma,\phi)$--conformal density, then for $x=p$, $\nu=\frac{\nu_p}{\norm{\nu_p}}$ is a  $(\Gamma,\phi)$--conformal measure. The converse also holds:

\begin{lemma}
Given a $(\Gamma,\phi)$--conformal measure $\nu$, there exists a unique conformal density $\set{\nu_x}$ such that $\nu_p=\nu$ for some choice of basepoint ${p}$.
\end{lemma}

\begin{proof}
Since $\nu_p=\nu$, then $\nu_x$ (if it exists) must be uniquely defined via
\[d\nu_x({[k]})=e^{-\phi(\beta_{[k]}(x,p))}d\nu({[k]})\]
for any $x\in X$ and ${[k]}\in \mc{F}$. Since $\beta_{[k]}(x,p)$ is bounded and continuous in $[k]$ for $x$ and $p$ fixed, and $\nu$ is probability, this defines a family of finite nonzero Borel measures $\{\nu_x\}$ and now we check that it satisfies the two remaining properties of a $(\Ga,\phi)$--conformal density:
\begin{enumerate}
\item For any $x,y\in X$ and ${[k]}\in \mc{F}$, we have
\[\frac{d\nu_x}{d\nu_y}({[k]})=\frac{e^{-\phi(\beta_{[k]}(x,p))}}{e^{-\phi(\beta_{[k]}(y,p))}}=e^{-\phi(\beta_{[k]}(x,p))-\beta_{[k]}(y,p))}=e^{-\phi(\beta_{[k]}(x,y))}.\]
\item For any $x\in X$, $\gamma\in \Ga$ and ${[k]}\in \mc{F}$, we have
\begin{align*}
\ga_*d\nu_x({[k]})=d\nu_x(\ga^{-1}{[k]})&=e^{-\phi(\beta_{\ga^{-1}{[k]}}(x,p))}d\nu(\ga^{-1}{[k]})=e^{-\phi(\beta_{[k]}(\ga x,\ga p))}(\ga_*d\nu)({[k]})\\
&=e^{-\phi(\beta_{[k]}(\ga x,p))}e^{\phi(\beta_{[k]}(\ga p, p))}(\ga_*d\nu)({[k]})=e^{-\phi(\beta_{[k]}(\ga x,p))}d\nu({[k]})\\
&=d\nu_{\ga x}({[k]}).
\end{align*}
\end{enumerate}
This proves the lemma.
\end{proof}

For $z\in \mf{a}^+$, and such that $\psi_\Gamma(z)\neq 0$, we say that $\phi \in \mathfrak{a}^{*}$ \textbf{extremizes $\psi_\Gamma$ at $z$,} or equivalently that it is \textbf{$\Gamma$--extremal at $z$}, if it is tangent to $\psi_{\Gamma}$ at $z$ in the sense that $\phi \geq \psi_{\Gamma}$ and $\phi(z)=\psi_{\Gamma}(z)$. If $z$ is interior to the non-negative support of $\psi_\Ga$, then extremal $\phi$ at $z$ always exist. The following theorem shows the existence of Patterson--Sullivan measures in higher rank.

\begin{thm}[{\cite[Thm 8.4]{Quint02}}]\label{thm:patt-sull-exists}
If $\phi$ extremizes $\psi_\Ga$ at a point $z\in \mf{a}^{++}$, then there exists a $(\Ga,\phi)$--density $\set{\nu_x}_{x\in X}$ supported on $\Lambda_{\mc{F}}(\Ga)\subset \mc{F}$.
\end{thm}

Since the above interpretation is purely Lie theoretic, we now relate it to the usual $\delta$--conformal density for Hadamard spaces. Given $z\in \mf{a}$ and $\xi_z=\lim_{t\to\infty}\exp(tz)\cdot p$, let $i_z\colon \mc{F}\to G\cdot \xi_z$ be the orbit map given by $[k]\mapsto k\xi_z$ (this is well defined since $M$ stabilizes $\xi_z$). When $z$ is a regular vector (i.e.~$z\in W\cdot\mf{a}^{++}$, where $W$ is the Weyl group), then $i_z$ is a homeomorphism. Otherwise, it is a bundle map whose fibre is homeomorphic to $M_z/M$ where $M_z$ is the stabilizer of $z$ in $K$.

\begin{prop}\label{prop:boundary_meas}
Let \(\phi\in\mathfrak a^*\setminus\{0\}\) extremize
\(\psi_\Gamma\) at \(z\in\mathfrak a^{++}\), and assume
$z'=\frac{\phi^*}{\norm{\phi^*}}\in {\mf{a}^+}$.
Let
\[
 i_{z'}\colon \mc{F}=G/P\longrightarrow G\cdot \xi_{z'}
\]
be the natural $G$--equivariant quotient.  If
\(\{\nu_x\}_{x\in X}\) is a \((\Gamma,\phi)\)-conformal density on
\(\mc{F}\), then $\mu_x=(i_{z'})_*\nu_x$ is a $\delta_\phi(\Gamma)$--conformal density on $G\cdot \xi_{z'}$, and $\supp(\mu_x)\subset i_{z'}(\Lambda_{\mc{F}}(\Gamma))$.


\end{prop}

\begin{proof}
Since \(\phi^*=\norm{\phi}\,z'\), one has
\[
 \phi(H)=\delta_\phi(\Gamma)\inner{z',H}
 \qquad(H\in\mathfrak a).
\]
Evaluating at \(z\) gives
\[
 \delta_\phi(\Gamma)
 =\frac{\phi(z)}{\inner{z',z}}
 =\frac{\psi_\Gamma(z)}{\inner{z',z}}.
\]
For \(z'\in{\mf{a}^+}\), the ordinary and
vector-valued Busemann cocycles satisfy
\begin{equation}\label{eq:busemann-pairing-audit}
 \inner{\beta_{[k]}(x,y),z'}
   =B_{i_{z'}([k])}(x,y).
\end{equation}

To prove it is a $\delta_\phi(\Ga)$--conformal density, we will use the definition of a $(\Gamma,\phi)$--conformal density. We first observe that the $G$--action commutes with $i_{z'}$ on pullbacks of subsets of $G\cdot \xi_{z'}\subset \pa_\infty X$ since $M_\theta/M_\Pi$ fibers are preserved. Thus for $\ga\in \Gamma$ we compute for $A\subset G\cdot \xi_{z'}$, 
\begin{align*}
\ga_*\mu_x(A)&=i_{z'}\nu_x(\ga^{-1}A)=\nu_x(i_{z'}^{-1}\ga^{-1}A)=\nu_x(\ga^{-1}i_{z'}^{-1}A)\\
&=\ga_*\nu_x(i_{z'}^{-1}A)=\nu_{\ga x}(i_{z'}^{-1}A)=\mu_{\ga x}(A).
\end{align*}

By conformality of \(\nu\),
\[
 \frac{d\nu_x}{d\nu_y}([k])
 =e^{-\phi(\beta_{[k]}(x,y))}
 =e^{-\delta_\phi(\Gamma) B_{i_{z'}([k])}(x,y)}
\]
for \(\nu_y\)--almost every \([k]\).  Therefore, for every Borel set
\(A\subset G\xi_{z'}\),
\begin{align*}
 \mu_x(A)
 &=\nu_x(i_{z'}^{-1}A)=\int_{i_{z'}^{-1}A}e^{-\delta_\phi(\Gamma) B_{i_{z'}([k])}(x,y)}\,d\nu_y([k])=\int_A e^{-\delta_\phi(\Gamma) B_\xi(x,y)}\,d\mu_y(\xi).
\end{align*}
Consequently,
\[
 \frac{d\mu_x}{d\mu_y}(\xi)
 =e^{-\delta_\phi(\Gamma) B_\xi(x,y)}
\]
for \(\mu_y\)--almost every \(\xi\).

For the conformality in the general case, we compute for $x,y\in X$ and $\xi\in G\xi_{z'}$,
\begin{align*}
\frac{d\mu_x}{d\mu_y}(\xi)&=\lim_{\eps\to 0} \frac{\mu_x(B_\eps(\xi))}{\mu_y(B_\eps(\xi))} =\lim_{\eps\to 0}\frac{\int_{i_{z'}^{-1}B_\eps(\xi)} e^{-\delta_\phi(\Gamma) B_{k\xi}(x,y)}d\nu_y([k])}{\nu_y(i_{z'}^{-1}B_\eps(\xi))}.
\end{align*}
As $k\xi$ is invariant under the conjugate of $M_\theta$ (recall $M_\theta$ fixes $w^{-1}\xi_{z'}$), the Busemann function at $k\xi$ only depends on the class $[k]\in K/M_\theta$. When $\eps$ is sufficiently small, the $G$--orbit of $\xi$ is close to its tangent space. Hence, by continuity of the Busemann function, for $k\xi\in B_\eps(\xi)$ we have $B_{k\xi}(x,y)=B_{\xi}(x,y)+O(\eps)$ where the constant may depend on $\xi_{z'}$. Therefore, as $\eps\to 0$ this limits to $\frac{d\mu_x}{d\mu_y}(\xi)=e^{-\delta_\phi(\Gamma) B_{\xi}(x,y)}$ as expected.
\end{proof}

\begin{remark}
For our applications we only need the case that $z'\in {\mf{a}^{++}}$. (The proof in this case reduces substantially.) Note that while we used the assumption $z'\in \mf{a}^+$, we suspect the statement holds for any $z'\in \mf{a}$ arising from $\phi^*$.
\end{remark}

Quint proved the following special case of Proposition \ref{prop:boundary_meas} when $z$ is the unique direction where the indicator function achieves the supremum. In this case, we can choose $z'$ in the same direction of $z$, and that $\delta_\phi(\Gamma)=\psi_\Gamma(\frac{z}{||z||})=\delta(\Gamma)$, the ordinary critical exponent. This in turn is a generalization of \cite[Thm A]{Alb99}.

\begin{prop}[{\cite[Prop 8.10]{Quint02}}]\label{prop:exist_conf}
Let $u\in \mc{C}_\Ga$ be the unique unit vector such that 
\[ \psi_\Ga(u)=\sup_{\substack{v\in \mc{C}_\Ga}{\norm{v}=1} } \psi_\Ga(v). \]
Then there exists at least one $\delta(\Ga)$--conformal density $\set{\mu_x}_{x\in X}$ supported on $\Lambda(\Ga) \cap G\cdot \xi_u\subset \pa_\infty X$.
\end{prop}

\begin{remark}
The support of a $\delta(\Gamma)$--conformal density arising from the above can lie on any $G$--orbit. It is shown in \cite{Benoist97} that given any closed cone $\Omega\subset\mf{a}^+$ that is invariant under a certain involution, there exists a discrete, Zariski-dense subgroup such that $\mc{C}_\Ga=\Omega$. Moreover, these groups can be taken to be Schottky (free) groups and hence have cohomological dimension $1$.
\end{remark}

\subsection{Applications}

By the above discussions, we can now give an explicit upper bound on $j_X(\Ga)$. The most general estimate uses Proposition \ref{prop:boundary_meas} where we have the freedom to choose arbitrary $z\in \mathcal C_\Ga$, and arbitrary $\phi\in \mathfrak a^*$ which extremizes $z$. For each of the choices, the Patterson--Sullivan measures vary both on the dimension of the conformal density and on the support of the $G$--orbit.

We consider the special realization of the Furstenberg boundary $\partial_F X=G\xi_z$ as the entire $G$--orbit of $\xi_z$, where $z\in \mathfrak a^{++}$ is to be specified such that it is extremized by some $\phi\in \mathfrak a^*$. Following Proposition \ref{prop:boundary_meas}, it determines a unique unit vector $z'\in \mathfrak a$, and there is a $\delta_{\phi}(\Ga)$--conformal density whose support is on the $G$--orbit of $\xi_{z'}$. 

Note that if $\psi_\Ga$ is $C^1$ at $z$ then there is a unique $\phi$ that extremizes $\psi_\Ga$ at $z$. At singular points of $\psi_\Ga$ there is a continuum of extremal $\phi$. For each pair $(z,\phi)$, Proposition \ref{prop:boundary_meas} shows that there exists a $\delta_\phi(\Ga)$--conformal density $\mu$ whose support lies in $G\cdot \xi_{z'}$. Recall the definition of  $j_X(\Gamma,\mu)\in\N$ from \cref{D:CritIndex}. 
Since $\partial_\mu X\subset G\cdot \xi_{z'}$, we know that the eigenvalues of $\nabla dB_{(x,\theta)}$ are independent of $x$ and $\theta$. To see this, we take $g\in G$ to be the isometry which sends the basepoint $p$ to $x$ and $z'$ to $v_{x\theta}$. Since the isometry preserves the second fundamental form of the horospheres, the Hessian of $\nabla dB_{(x,\theta)}$ and $\nabla dB_{(p,z')}$ have the same eigenvalues.

For higher rank symmetric spaces, given a direction $u\in \mathfrak a\subset \mathfrak p$, it is well known that 
\[\nabla dB_{(p,u)}=\sqrt{-R_u}=\sqrt{(\ad_u)^2\rest{\mathfrak{p}}},\]
where $R_u=R(u,\cdot)u$ is the Jacobi operator, and $\mathfrak g=\mathfrak k+\mf{p}$ is the Cartan decomposition so that $\mathfrak p$ identifies with the tangent space $T_pX$. Note that the root space decomposition $\mathfrak p=\mathfrak a+\bigoplus_{\alpha\in \Lambda^+}\mathfrak p_\alpha$ diagonalizes $R_u$ hence also does $\nabla dB_{(p,u)}$. Thus, $\nabla dB_{(p,u)}$ has at least $r$ zero eigenvalues corresponding to $\mathfrak a$ and all other eigenvalues are nonnegative, listed as $\{|\alpha_i(u)|:\alpha_i\in \Lambda^+\}$, where $\Lambda^+$ is any choice of positive roots, and the roots are counted with multiplicity. Thus, if we choose $\Lambda^+$ to be a positive system with respect to $u$ and denote $L_u$ the diagonal matrix given by
\[L_u = \op{diag}(0^{(r)},\alpha_1(u),...,\alpha_{n-r}(u)),\]
where $n,r$ are the dimension and rank of the symmetric space, then our interest is in
\[ \inf_{\begin{subarray}{l} \quad\;\;\; z\in \mathcal C_\Ga\\
\phi \text{ extremizes } z
\end{subarray}}\min\{k:\op{tr}_k L_{z'}-\delta_\phi(\Ga)>0\}. \] 
This gives an upper bound on $j_X(\Ga)$ and hence provides our optimal bound for the homological dimension. We summarize the discussions into the following theorem.

\begin{thm}\label{thm:j_X}
Let $X=G/K$ be a higher rank symmetric space of non-compact type, and $\Ga<G$ be a discrete, Zariski-dense subgroup. Then
\[j_X(\Ga)\leq  \inf_{\begin{subarray}{l} \quad\;\;\; z\in \mathcal C_\Ga\cap\mathfrak{a}^{++}\\ \phi \text{ extremizes } z \end{subarray}}\min\left(\{k:\op{tr}_k L_{z'}-\delta_\phi(\Ga)>0\}\cup\{n+1\}\right).\]
\end{thm}

\subsection{Explicit bounds} 

We now let $G$ be a real, simple Lie group of real rank at least two, and $\Ga<G$ be an infinite covolume, Zariski-dense subgroup. Since $G$ has Property (T), there is a uniform gap on the critical exponent of $\Ga$ (see \cite{Leuzinger03}), that is, there exists $\eps_0(G)>0$ such that $\delta(\Ga)<h-\eps_0(G)$ where $h$ is the volume entropy of the corresponding symmetric space $X$ (which coincides with the critical exponent of a lattice). Moreover, Quint \cite{Quint03} proved there is a gap in terms of the indicator function, and an explicit bound $\psi_\Ga\leq 2\rho-\Theta$ was given by Lee--Oh \cite[Thm 7.1]{LeeOh22}, where $\rho$ is the half sum of all positive roots (counted with multiplicity) of $(\mf{g},\mf{a}^+)$ and $\Theta$ is the half sum of roots in a maximal strongly orthogonal system of $(\mf{g},\mf{a}^+)$. For a given choice of $\mf{a}^+$, $\Theta$ is independent of the choice of maximal strongly orthogonal system. The maximal orthogonal systems are explicitly computed in \cite[App]{Oh02}. 

The Anosov case is closely tied to the entropy, pressure, spectral, and growth-indicator theory for higher rank representations
\cite{Sambarino14,Sambarino14Counting,
BridgemanCanaryLabourieSambarino15,KimEtAl21,
PozzettiSambarinoWienhard23,LutskoWeichWolf24,Wolf25LimitCone}.
Kim--Minsky--Oh in \cite{KimEtAl21} conjectured that if $\Gamma$ is Borel--Anosov, then $\psi_\Gamma\leq \rho$. This conjecture has now been proved in \cite{LutskoWeichWolf24} except for the cases $\mathfrak{sl}_3(\mathbf{K})$, $\mathbf{K}=\mathbf{R},\mathbf{C},\mathbf{H}$, and $\mathfrak e_6^{-26}$. More generally, Wolf in \cite{Wolf25LimitCone} proved that $\psi_\Gamma\leq\rho$ whenever the limit cone of $\Gamma$ avoids two Weyl chamber facets which are not interchanged by the opposition involution. In particular, this applies to every $I$--Anosov subgroup for which $I$ contains at least two distinct simple roots not interchanged by the opposition involution. If $\Gamma<G$ is Borel--Anosov, then its Cartan projections are uniformly regular, and its limit cone lies in the interior of the positive Weyl chamber. Thus the regularity assumptions in the growth-indicator estimates below are automatic for such groups.


The relationship between Anosov dynamics, entropy, counting, and
higher rank Patterson--Sullivan theory has been developed in
\cite{GueritaudGuichardKasselWienhard17,GuichardKasselWienhard18,BridgemanCanaryLabourieSambarino15,
EdwardsOh22,EdwardsLeeOh20,EdwardsLeeOh21,LeeOh22}. More recently, Wolf and his coauthors have clarified the analytic significance of the growth-indicator function. Wolf--Zhang relate Quint's growth indicator to the modified critical exponent for the polyhedral distance and use this to characterize the bottom of the $L^2$--spectrum of the Laplace--Beltrami operator in higher rank \cite{WolfZhang23}. Building on this perspective, Lutsko--Weich--Wolf prove a precise relation between growth rates of $\Gamma$, polyhedral bounds on the joint spectrum of invariant differential operators, and decay of matrix coefficients. In particular, they obtain a general characterization of temperedness of $L^2(\Gamma\backslash G)$ \cite{LutskoWeichWolf24}. Thus the growth-indicator function $\psi_\Gamma$ provides a bridge
between the asymptotic geometry of Cartan projections and spectral or representation-theoretic properties of locally symmetric spaces. The estimates used in Theorem \ref{thm:anosov} exploit precisely this bridge
to translate higher rank asymptotic geometry into the homological restrictions proved below.

The estimates below therefore interpolate between the general $\Theta$--gap of Lee--Oh and the stronger $\rho$--gap phenomena expected, and in many cases known, for Anosov subgroups. This places the critical index estimate in the same circle of ideas as the higher rank entropy and growth-indicator theory of \cite{Sambarino14,Sambarino14Counting,BridgemanCanaryLabourieSambarino15,KimEtAl21,LeeOh22}.

In what follows, we assume that $\psi$ has an $\eta$--gap (i.e., $\psi_\Gamma\leq 2\rho-\eta$). The following estimate was partly suggested by \cite[Prop 5.27]{Fraczyk-Lowe}, where they obtained the result for the case of $G=\SL(r,\mathbf{R})$ and $\eta=\Theta$.

\begin{thm}\label{thm:explicit upper bound}
Let $G$ be a simple Lie group with real rank at least two and $\Ga$ be a discrete, torsion-free, Zariski-dense subgroup. Suppose
\begin{enumerate}
\item $\psi_\Gamma\leq 2\rho-\eta$ on $\mathfrak a^+$, where $\eta\in \mathfrak a^*$ satisfies $(2\rho-\eta)^*\in \mathfrak a^+$, and
\item there exists $\lambda\in (0,1]$ such that $\lambda(2\rho-\eta)$ extremizes $\psi_\Gamma$ at a point $z \in\mathfrak a^{++}$,
\end{enumerate}
then 
\begin{equation}\label{eq:j_X_bound}
j_X(\Ga)\leq n+1-\frac{\inner{\eta,2\rho-\eta}}{\ell}
\end{equation}
where 
\[\ell=\max_{\alpha\in \Sigma^+}\inner{\alpha,2\rho-\eta}.\] 
\end{thm}

\begin{remark}
Note that if $\mathcal C_\Gamma\subset \mathfrak a^{++}$ (for example if $\Ga$ is Borel--Anosov), then the second condition of the above theorem is automatically satisfied. We needed this second condition, ensuring the extremal point belongs to $\mathfrak a^{++}$, in order to satisfy the condition in \cref{thm:patt-sull-exists} for the existence of $(\Gamma,\phi)$--densities. We are not certain whether this condition is in fact needed.
\end{remark}

\begin{proof}
Since $j_X(\Gamma)\leq n+1$, we may assume that $\langle \eta,\rho-\eta\rangle>0$. We set $\phi=\lambda(2\rho-\eta)$ and $z'=\phi^*/{|\phi^*|}=(2\rho-\eta)^*/|2\rho-\eta|\in \mathfrak a^+$. Now we apply Theorem \ref{thm:j_X} and get
\[j_X(\Ga)\leq \min\{k:\op{tr}_k L_{z'}>\delta_\phi(\Ga)\},\]
where from Proposition \ref{prop:boundary_meas} we have
\[\delta_\phi(\Gamma)=\frac{\psi_\Ga(z)}{\langle z',z\rangle}\leq \frac{(2\rho-\eta)(z)}{\frac{2\rho-\eta}{|2\rho-\eta|}(z)}=|2\rho-\eta|.\]
Thus, if we set $\{\alpha_i|r+1\leq i\leq n\}$ the set of all positive roots in the order
\[ \alpha_{r+1}(z')\leq \alpha_{r+2}(z')\leq\dots \leq \alpha_{n}(z'), \] 
and $k$ is some integer satisfying the following inequality
\begin{equation}\label{eq:pf-key}
\alpha_{r+1}(z')+\dots+\alpha_k(z')>|2\rho-\eta|,
\end{equation}
then we will have $j_X(\Ga)\leq k$. Note that $\alpha_{r+1}+\dots+\alpha_{n}=2\rho$, so the inequality \eqref{eq:pf-key} is equivalent to
\begin{equation}\label{pf:main-ineq}
\langle\alpha_{k+1}+\dots+\alpha_{n}, 2\rho-\eta\rangle<\langle \eta, 2\rho-\eta\rangle.
\end{equation}
Since $\ell=\max_{\alpha\in \Sigma^+}\inner{\alpha,2\rho-\eta}$, if we choose
\[k=\begin{cases} n+1-\frac{\langle \eta, 2\rho-\eta\rangle}{\ell} & \textrm{ if }\frac{\langle \eta, 2\rho-\eta\rangle}{\ell}\in \mathbf{N}\\ n-\floor{\frac{\langle \eta, 2\rho-\eta\rangle}{\ell}} & \textrm{ if }\frac{\langle \eta, 2\rho-\eta\rangle}{\ell}\notin \mathbf{N} \end{cases},\]
then the above inequality \eqref{pf:main-ineq} holds. Hence, we have
\begin{equation*}
j_X(\Gamma)\leq k\leq n+1-\frac{\langle \eta, 2\rho-\eta\rangle}{\ell}.
\end{equation*}
\end{proof}

In the following subsection, we present two explicit bounds (for $\eta=\Theta$ and $\eta=\rho$) of $j_X(\Ga)$ in relation to certain strengthened conditions imposed on $\psi$.

\subsubsection{\texorpdfstring{$\Theta$}{\Theta}--gap and \texorpdfstring{$\rho$}{rho}--gap bounds}
We now apply the above theorem to the cases $\eta=\Theta$ and $\eta=\rho$. This leads to the following theorems. We defer the proofs to Section \ref{sec:proof_general}, \ref{sec:anosov} and Appendix \ref{app:allcases}.

\begin{thm}\label{thm:general}
Let $G$ be a non-compact, real, simple Lie group with real rank $r$, let $n$ be the dimension of the associated symmetric space $G/K$, and let $\Gamma<G$ be a discrete, finitely generated, torsion-free, Zariski-dense subgroup which is not a lattice. Suppose further there exists $\lambda\in (0,1]$ such that $\lambda(2\rho-\Theta)$ extremizes $\psi_\Gamma$ at a point $z \in\mathfrak a^{++}$. Then we have the following upper bound on $j_X(\Gamma)$ according to the type of restricted root system of $G$.  
\begin{enumerate}
\item If $G$ is of Type $\mathrm{A}_r$, then $j_X(\Gamma)< n+1-\frac{r}{8}$.
\item If $G$ is of Type $\mathrm{B}_r$ or $D_r$, then $j_X(\Gamma)< n+1-\frac{3r}{16}$.
\item If $G$ is of Type $\mathrm{C}_r$ or $(\mathrm{BC})_r$, then $j_X(\Gamma)< n+1-\frac{r}{4}$.
\end{enumerate}
In all classical cases of $G$, we have $j_X(\Gamma)< n+1-\frac{r}{8}$.
\end{thm}

\begin{thm}
Let $G$ be a non-compact, real, simple Lie group with real rank at least two and which is not of type $\mathrm{A}_2$ restricted root system, $n$ be the dimension of the associated symmetric space $G/K$, and $\Gamma<G$ be a discrete, finitely generated, torsion-free, Zariski-dense subgroup with regular limit cone (i.e. $\mathcal C_\Ga\subset \mathfrak a^{++}\cup\set{0}$). Then we have  
\[ j_X(\Gamma)< \frac{5n}{6}+1. \]
In particular, if $\Ga$ is a Borel--Anosov subgroup of $G$, then this holds.
\end{thm}

\section{Proof of Theorem \ref{thm:j_X_bound}: Real split cases}\label{sec:proof_general}

Since $\Ga$ is not a lattice, we have $\psi_\Ga\leq 2\rho-\Theta$ by \cite[Thm 7.1]{LeeOh22}. We now apply Theorem \ref{thm:explicit upper bound} for the case $\eta=\Theta$. We first note that the quantity $\frac{\langle \Theta, 2\rho-\Theta\rangle}{\ell}$ is positive and since $j_X(\Ga)$ is an integer, we have $j_X(\Ga)\leq n$ in all cases. This indeed proves the theorem for low ranks $r=2$ and $r=3$.

For the remaining cases, we assume $r\geq 4$, and it suffices to compute case by case the quantity $\frac{\langle \Theta, 2\rho-\Theta\rangle}{\ell}$. We will also check along the way that $(2\rho-\Theta)^*\in \mathfrak a^+$.

For readability, here we only present the proof for the case when $G$ is real split. These are exactly the cases where all roots have multiplicity one. However, the computations are representable enough, and in fact the final conclusions are the same within each type of Dynkin diagram. We simply leave the remaining cases to the Appendix \ref{app:allcases}.

According to \cite[App C]{Knapp} and \cite[App]{Oh02}, the root data (including the set of simple roots $\{\alpha_1,\dots, \alpha_r\}$, the set of positive roots $\Delta^+$, the sum of all positive roots $2\rho$), and the half sum of the roots in a maximal strongly orthogonal system $\Theta$ is listed in Table \ref{tab:root_data}.

\begin{table}
\centering
\begin{tabular}{|c|c|c|}
           \hline
            $(\mathfrak g, \mathfrak a)$ & $\mathrm{A}_r$ & $\mathrm{B}_r$ \\
            \hline
            $\{\alpha_1,\dots, \alpha_r\}$& $\{e_1-e_2,e_2-e_3, \dots, e_r-e_{r+1}\}$ & $\{e_1-e_2,e_2-e_3, \dots, e_{r-1}-e_r, e_r\}$ \\
            \hline
            $\Delta^+$ & $\{e_i-e_j|i<j\}$ & $\{e_i\pm e_j|i<j\}\cup \{e_i\}$ \\
            \hline
            $2\rho$ & $re_1+(r-2)e_2+\dots+(-r)e_{r+1}$  & $(2r-1)e_1+(2r-3)e_2+\dots+e_{r}$  \\
            \hline
            $\Theta$ if $r$ is odd & $\sum_{i=1}^{(r-1)/2}\frac{i}{2}\alpha_i+\sum_{i=(r+1)/2}^{r}\frac{r-i+1}{2}\alpha_i$ 
            & $\sum_{i=1}^{(r-1)/2}i\alpha_i+\sum_{i=(r+1)/2}^{r}\frac{r}{2}\alpha_i$  \\
            \hline
            $\Theta$ if $r$ is even & $\sum_{i=1}^{r/2}\frac{i}{2}\alpha_i+\frac{r}{4}\alpha_{r/2+1}+\sum_{i=r/2+2}^{r}\frac{r-i+1}{2}\alpha_i$ 
            & $\sum_{i=1}^{r/2}i\alpha_i+\sum_{i=r/2+1}^{r}\frac{r}{2}\alpha_i$ \\
            \hline
\end{tabular}

\quad

\begin{tabular}{|c|c|c|}
           \hline
            $(\mathfrak g, \mathfrak a)$ & $\mathrm{C}_r$ & $\mathrm{D}_r$ \\
            \hline
            $\{\alpha_1,\dots, \alpha_r\}$& $\{e_1-e_2,e_2-e_3, \dots, e_{r-1}-e_r, 2e_r\}$ & $\{e_1-e_2,e_2-e_3, \dots, e_{r-1}-e_r, e_{r-1}+e_r\}$ \\
            \hline
            $\Delta^+$ & $\{e_i\pm e_j|i<j\}\cup \{2e_i\}$ & $\{e_i\pm e_j|i<j\}$ \\
            \hline
            $2\rho$ & $(2r)e_1+(2r-2)e_2+\dots+2e_{r}$  & $(2r-2)e_1+(2r-4)e_2+\dots+2e_{r-1}$  \\
            \hline
            $\Theta$ if $r$ is odd & $\sum_{i=1}^{r-1}i\alpha_i+\frac{r}{2}\alpha_{r}$
            & $\sum_{i=1}^{\frac{r-1}2}i\alpha_i+\sum_{i=\frac{r+1}2}^{r-2}\frac{r-1}{2}\alpha_i+\frac{r-1}{4}(\alpha_{r-1}+\alpha_r)$  \\
            \hline
            $\Theta$ if $r$ is even & $\sum_{i=1}^{r-1}i\alpha_i+\frac{r}{2}\alpha_{r}$ 
            & $\sum_{i=1}^{r/2}i\alpha_i+\sum_{i=r/2+1}^{r-2}\frac{r}{2}\alpha_i+\frac{r}{4}(\alpha_{r-1}+\alpha_r)$ \\
            \hline
\end{tabular}\caption{Root data of $G$}\label{tab:root_data}
\end{table}

\begin{enumerate}
\item {Type $\mathrm{A}_r$}: the real split Lie group is $G=\SL(r+1,\mathbf{R})$. This is already proved in \cite[Prop 5.27]{Fraczyk-Lowe}. To fit in our language and also for the completeness, we include it here. When $r=2m+1$ is odd, we compute from Table \ref{tab:root_data} that 
\[ \Theta=\frac{1}{2}e_1+\dots+\frac{1}{2}e_{m+1}-\frac{1}{2}e_{m+2}-\dots-\frac{1}{2}e_{r+1} \]
and
\begin{align*}
2\rho-\Theta&=\frac{4m+1}{2}e_1+\frac{4m-3}{2}e_2+\dots+\frac{1}{2}e_{m+1}-\frac{1}{2}e_{m+2}-\dots-\frac{4m-3}{2}e_r-\frac{4m+1}{2}e_{r+1}\\ &\in \mathfrak a^+.
\end{align*}
Thus we have
\begin{align*}
\langle \Theta, 2\rho-\Theta\rangle&=\left(\frac{4m+1}{4}+\frac{4m-3}{4}+\dots+\frac{1}{4}\right)\times 2=\frac{r(r+1)}{4},
\end{align*}
and $\ell=\langle e_1-e_{r+1}, 2\rho-\Theta\rangle=2r-1$. Therefore, by \eqref{eq:j_X_bound},
\[j_X(\Gamma)\leq n+1-\frac{\langle \Theta, 2\rho-\Theta\rangle}{\ell} < n+1-\frac{r}{8}.\]
Similarly when $r=2m$ is even, we have
\[ \Theta=\frac{1}{2}e_1+\dots+\frac{1}{2}e_{m}+0\cdot e_{m+1}-\frac{1}{2}e_{m+2}-\dots-\frac{1}{2}e_{r+1} \]
and
\begin{align*}
2\rho-\Theta&=\frac{4m-1}{2}e_1+\frac{4m-5}{2}e_2+\dots+\frac{3}{2}e_{m}+0\cdot e_{m+1}-\frac{3}{2}e_{m+2}-\dots-\frac{4m-5}{2}e_r-\frac{4m-1}{2}e_{r+1}\\ &\in \mathfrak a^+.
\end{align*}
Thus we have
\begin{align*}
\langle \Theta, 2\rho-\Theta\rangle&=\left(\frac{4m-1}{4}+\frac{4m-5}{4}+\dots+\frac{3}{4}\right)\times 2=\frac{r(r+1)}{4},
\end{align*}
and $\ell=\langle e_1-e_{r+1}, 2\rho-\Theta\rangle=2r-1$. Therefore, by \eqref{eq:j_X_bound},
\[j_X(\Gamma)\leq n+1-\frac{\langle \Theta, 2\rho-\Theta\rangle}{\ell} < n+1-\frac{r}{8}.\]

\item Type $\mathrm{B}_r$: the real split Lie group is $G=\SO(r,r+1)$. When $r=2m+1$ is odd, we compute from Table \ref{tab:root_data} that 
\[ \Theta=e_1+\dots+e_{m}+\frac{1}{2}e_{m+1} \] 
and
\begin{align*}
2\rho-\Theta&=(4m)e_1+(4m-2)e_2+\dots+(2m+2)e_{m}+(2m+\frac{1}{2})e_{m+1}+(2m-1)e_{m+2}+\dots+e_{2m+1}\\ &\in \mathfrak a^+.
\end{align*}
Thus we have
\begin{align*}
\langle \Theta, 2\rho-\Theta\rangle&=\left(4m+(4m-2)+\dots+(2m+2)\right)+(m+\frac{1}{4})=\frac{r(3r-2)}{4},
\end{align*}
and $\ell=\langle e_1+e_2, 2\rho-\Theta\rangle=4r-6$. Therefore, by \eqref{eq:j_X_bound},
\[j_X(\Gamma)\leq n+1-\frac{\langle \Theta, 2\rho-\Theta\rangle}{\ell} < n+1-\frac{3r}{16}.\]
When $r=2m$ is even, we compute from Table \ref{tab:root_data} that 
\[ \Theta=e_1+\dots+e_{m} \] 
and
\begin{align*}
2\rho-\Theta&=(2r-2)e_1+(2r-4)e_2+\dots+(r)e_{m}+(r-1)e_{m+1}+(r-3)e_{m+2}+\dots+e_{r}\\
&\in \mathfrak a^+.
\end{align*}
Thus we have
\begin{align*}
\langle \Theta, 2\rho-\Theta\rangle&=(2r-2)+(2r-4)+\dots+r=\frac{r(3r-2)}{4},
\end{align*}
and $\ell=\langle e_1+e_2, 2\rho-\Theta\rangle=4r-6$. Therefore, by \eqref{eq:j_X_bound},
\[j_X(\Gamma)\leq n+1-\frac{\langle \Theta, 2\rho-\Theta\rangle}{\ell} < n+1-\frac{3r}{16}.\]

\item Type $\mathrm{C}_r$: the real split Lie group is $G=\Sp(r,\mathbf{R})$. We compute from Table \ref{tab:root_data} that 
\[ \Theta=e_1+\dots+e_r \]
and
\[2\rho-\Theta=(2r-1)e_1+(2r-3)e_2+\dots+e_r\in \mathfrak a^+.\]
Thus we have
\begin{align*}
\langle \Theta, 2\rho-\Theta\rangle&=(2r-1)+(2r-3)+\dots+1=r^2,
\end{align*}
and $\ell=\langle 2e_1, 2\rho-\Theta\rangle=4r-2$. Therefore, by \eqref{eq:j_X_bound},
\[j_X(\Gamma)\leq n+1-\frac{\langle \Theta, 2\rho-\Theta\rangle}{\ell} < n+1-\frac{r}{4}.\]
\item Type $\mathrm{D}_r$: the real split Lie group is $G=\SO(r,r)$. When $r=2m+1$ is odd, we compute from Table \ref{tab:root_data} that 
\[ \Theta=e_1+\dots+e_m \]
and
\begin{align*}
2\rho-\Theta&=(2r-3)e_1+(2r-5)e_2+\dots+(r)e_{m}+(r-1)e_{m+1}+(r-3)e_{m+2}+\dots+2e_{r-1}\\ &\in \mathfrak a^+.
\end{align*}
Thus we have
\begin{align*}
\langle \Theta, 2\rho-\Theta\rangle&=(2r-3)+(2r-5)+\dots+r=\frac{3(r-1)^2}{4},
\end{align*}
and $\ell=\langle e_1+e_2, 2\rho-\Theta\rangle=4r-8$. Therefore, by \eqref{eq:j_X_bound},
\[j_X(\Gamma)\leq n+1-\frac{\langle \Theta, 2\rho-\Theta\rangle}{\ell} < n+1-\frac{3r}{16}.\]
    
When $r=2m$ is even, we compute from Table \ref{tab:root_data} that 
\[ \Theta=e_1+\dots+e_m \]
and
\begin{align*}
2\rho-\Theta&=(2r-3)e_1+(2r-5)e_2+\dots+(r-1)e_{m}+(r-2)e_{m+1}+(r-4)e_{m+2}+\dots+2e_{r-1}\\ &\in \mathfrak a^+.
\end{align*}
Thus we have
\begin{align*}
\langle \Theta, 2\rho-\Theta\rangle&=(2r-3)+(2r-5)+\dots+(r-1)=\frac{r(3r-4)}{4},
\end{align*}
and $\ell=\langle e_1+e_2, 2\rho-\Theta\rangle=4r-8$. Therefore, by \eqref{eq:j_X_bound},
\[j_X(\Gamma)\leq n+1-\frac{\langle \Theta, 2\rho-\Theta\rangle}{\ell} < n+1-\frac{3r}{16}.\]
\end{enumerate}

For the remaining cases in each type $\mathrm{A}_r-\mathrm{D}_r$ and type $(\mathrm{BC})_r$, the same inequality holds, and we leave the computations in Appendix \ref{app:allcases}. 

\section{Proof of Theorem \ref{thm:anosov}: Real split cases}\label{sec:anosov}

Since $\Ga$ has a regular limit cone, the second condition of Theorem \ref{thm:explicit upper bound} holds automatically. Moreover, since $G$ is not of type $\mathrm{A}_2$, we know that $\psi_\Ga\leq \rho$ by \cite{LutskoWeichWolf24}. Now we apply Theorem \ref{thm:explicit upper bound} for the case $\eta=\rho$, then we have
\[j_X(\Ga)\leq n+1-\frac{|\rho|^2}{\ell},\]
where $\ell=\max_{\alpha\in \Sigma^+}\langle \alpha,\rho\rangle$. Thus, it suffices to compute and estimate in all cases the quantity $\frac{|\rho|^2}{\ell}$. Similar to Section \ref{sec:proof_general}, we only present the proof for real split Lie groups and leave the remaining cases in Appendix \ref{app:allcases}. 

In the following Table \ref{tab:root_data_real_split}, we list the root data (including the set of simple roots $\{\alpha_1,\dots, \alpha_r\}$, the set of positive roots $\Delta^+$, the sum of all positive roots $2\rho$, and the dimension of the corresponding symmetric space $n$) for each type of real split Lie groups.
    
\begin{table}
\centering
\begin{tabular}{|c|c|c|}
            \hline
            $G$ & $\SL(r+1,\mathbf{R})$ & $\SO(r,r+1)$\\
           \hline
            $(\mathfrak g, \mathfrak a)$ & $\mathrm{A}_r$ & $\mathrm{B}_r$ \\
            \hline
            $\{\alpha_1,\dots, \alpha_r\}$& $\{e_1-e_2,e_2-e_3, \dots, e_r-e_{r+1}\}$ & $\{e_1-e_2,e_2-e_3, \dots, e_{r-1}-e_r, e_r\}$ \\
            \hline
            $\Delta^+$ & $\{e_i-e_j|i<j\}$ & $\{e_i\pm e_j|i<j\}\cup \{e_i\}$ \\
            \hline
            $2\rho$ & $re_1+(r-2)e_2+\dots+(-r)e_{r+1}$  & $(2r-1)e_1+(2r-3)e_2+\dots+e_{r}$  \\
            \hline
            $n$ & $\frac{r(r+3)}{2}$ & $r(r+1)$ \\
            \hline   
\end{tabular}

\quad

\begin{tabular}{|c|c|c|}
           \hline
            $G$ & $\Sp(r,\mathbf{R})$ & $\SO(r,r)$\\
           \hline
            $(\mathfrak g, \mathfrak a)$ & $\mathrm{C}_r$ & $\mathrm{D}_r$ \\
            \hline
            $\{\alpha_1,\dots, \alpha_r\}$& $\{e_1-e_2,e_2-e_3, \dots, e_{r-1}-e_r, 2e_r\}$ & $\{e_1-e_2,e_2-e_3, \dots, e_{r-1}-e_r, e_{r-1}+e_r\}$ \\
            \hline
            $\Delta^+$ & $\{e_i\pm e_j|i<j\}\cup \{2e_i\}$ & $\{e_i\pm e_j|i<j\}$ \\
            \hline
            $2\rho$ & $(2r)e_1+(2r-2)e_2+\dots+2e_{r}$  & $(2r-2)e_1+(2r-4)e_2+\dots+2e_{r-1}$  \\
            \hline
            $n$ & $r(r+1)$ & $r^2$ \\
            \hline
\end{tabular}

\quad

\begin{tabular}{|c|c|}
           \hline
            $(\mathfrak g, \mathfrak a)$ & $\mathrm{E}_6$  \\
            \hline
            $\{\alpha_1,\dots, \alpha_6\}$& $\{\frac{1}{2}(e_8-e_7-e_6-e_5-e_4-e_3-e_2+e_1), e_2+e_1, e_2-e_1, e_3-e_2,e_4-e_3,e_5-e_4\}$\\
            \hline
            $\Delta^+$ & $\{e_i\pm e_j|5\geq i>j\}\cup \{\frac{1}{2}(e_8-e_7-e_6+\sum_{i=1}^5(-1)^{n(i)}e_i)|\sum_i n(i) \textrm{ is even}\}$ \\
            \hline
            $2\rho$ & $2e_2+4e_3+6e_4+8e_5-8e_6-8e_7+8e_8$\\
            \hline
            $n$ & $42$\\
            \hline
\end{tabular}

\quad

\begin{tabular}{|c|c|}
           \hline
            $(\mathfrak g, \mathfrak a)$ & $\mathrm{E}_7$  \\
            \hline
            $\{\alpha_1,\dots, \alpha_7\}$& $\{\frac{1}{2}(e_8-e_7-e_6-e_5-e_4-e_3-e_2+e_1), e_2+e_1\}\cup\{e_{i+1}-e_i|1\leq i\leq 5\}\}$\\
            \hline
            $\Delta^+$ & $\{e_i\pm e_j|6\geq i>j\}\cup\{e_8-e_7\}\cup \{\frac{1}{2}(e_8-e_7+\sum_{i=1}^6(-1)^{n(i)}e_i)|\sum_i n(i) \textrm{ is odd}\}$ \\
            \hline
            $2\rho$ & $2e_2+4e_3+6e_4+8e_5+10e_6-17e_7+17e_8$\\
            \hline
             $n$ & $70$\\
            \hline
\end{tabular}

\quad

\begin{tabular}{|c|c|}
           \hline
            $(\mathfrak g, \mathfrak a)$ & $\mathrm{E}_8$  \\
            \hline
            $\{\alpha_1,\dots, \alpha_8\}$& $\{\frac{1}{2}(e_8-e_7-e_6-e_5-e_4-e_3-e_2+e_1), e_2+e_1\}\cup\{e_{i+1}-e_i|1\leq i\leq 6\}$\\
            \hline
            $\Delta^+$ & $\{e_i\pm e_j|i>j\}\cup \{\frac{1}{2}(e_8+\sum_{i=1}^7(-1)^{n(i)}e_i)|\sum_i n(i) \textrm{ is even}\}$ \\
            \hline
            $2\rho$ & $2e_2+4e_3+6e_4+8e_5+10e_6+12e_7+46e_8$\\
            \hline
             $n$ & $128$\\
            \hline
\end{tabular}

\quad 

\begin{tabular}{|c|c|}
           \hline
            $(\mathfrak g, \mathfrak a)$ & $\mathrm{F}_4$  \\
            \hline
            $\{\alpha_1,\dots, \alpha_4\}$& $\{\frac{1}{2}(e_1-e_2-e_3-e_4), e_4, e_3-e_4, e_2-e_3\}$\\
            \hline
            $\Delta^+$ & $\{e_i\pm e_j|i<j\}\cup \{e_i\}\cup \{\frac{1}{2}(e_1\pm e_2\pm e_3\pm e_4)\}$ \\
            \hline
            $2\rho$ & $11e_1+5e_2+3e_3+e_4$\\
            \hline
             $n$ & $28$\\
            \hline
\end{tabular}

\quad

\begin{tabular}{|c|c|}
           \hline
            $(\mathfrak g, \mathfrak a)$ & $\mathrm{G}_2$  \\
            \hline
            $\{\alpha_1,\alpha_2\}$& $\{e_1-e_2, -2e_1+e_2+e_3\}$\\
            \hline
            $\Delta^+$ & $\{e_1-e_2, e_3-e_1, e_3-e_2,-2e_1+e_2+e_3, e_1-2e_2+e_3, -e_1-e_2+2e_3\}$ \\
            \hline
            $2\rho$ & $-2e_1-4e_2+6e_3$\\
            \hline
             $n$ & $8$\\
            \hline
\end{tabular}\caption{Root data of $G$}\label{tab:root_data_real_split}
\end{table}

\begin{enumerate}
\item {Type $\mathrm{A}_r$}: For the real split Lie group $G=\SL(r+1,\mathbf{R})$, we see from Table \ref{tab:root_data_real_split} that 
\[ \rho=\frac{r}2e_1+\frac{r-2}2 e_2+\dots+(-\frac{r}2) e_{r+1}, \]
and that $\ell=\langle e_1-e_{r+1}, \rho\rangle=r$. Thus, when $r=2m+1$ is odd, we have
\begin{align*}
\frac{|\rho|^2}{\ell}&=\frac{\left[(\frac{2m+1}2)^2+(\frac{2m-1}{2})^2+\cdots+(\frac{1}{2})^2\right]\cdot 2}{r}=\frac{2m^2+5m+3}{6} >\frac{(2m+1)(m+2)}{6}=\frac{n}6.
\end{align*}
when $r=2m$ is even, we have
\begin{align*}
\frac{|\rho|^2}{\ell}&=\frac{\left[m^2+(m-1)^2+\cdots+1^2\right]\cdot 2}{2m}=\frac{2m^2+3m+1}{6} >\frac{2m^2+3m}{6}=\frac{n}6.  
\end{align*}

\item {Type $\mathrm{B}_r$}: For the real split Lie group $G=\SO(r,r+1)$, we see from Table \ref{tab:root_data_real_split} that 
\[ \rho=\frac{2r-1}2e_1+\frac{2r-3}2 e_2+\dots+\frac{1}2 e_{r}, \]
and that $\ell=\langle e_1+e_2, \rho\rangle=2r-2$. Thus, we have
\begin{align*}
\frac{|\rho|^2}{\ell}&=\frac{(\frac{2r-1}2)^2+(\frac{2r-3}{2})^2+\cdots+(\frac{1}{2})^2}{2r-2}=\frac{r(4r^2-1)}{24(r-1)} >\frac{r(r+1)}{6}=\frac{n}6.
\end{align*}

\item {Type $\mathrm{C}_r$}: For the real split Lie group $G=\Sp(r,\mathbf{R})$, we see from Table \ref{tab:root_data_real_split} that 
\[ \rho=r\cdot e_1+(r-1)\cdot e_2+\dots+1\cdot  e_{r}, \]
and that $\ell=\langle 2e_1, \rho\rangle=2r$. Thus, we have
\begin{align*}
\frac{|\rho|^2}{\ell}&=\frac{r^2+(r-1)^2+\cdots+1^2}{2r}=\frac{2r^2+3r+2}{12} >\frac{2r^2+2r}{12}=\frac{n}6.
\end{align*}
        
\item {Type $\mathrm{D}_r$}: For the real split Lie group $G=\SO(r,r)$, we see from Table \ref{tab:root_data_real_split} that 
\[ \rho=(r-1)\cdot e_1+(r-2)\cdot e_2+\dots+1\cdot  e_{r-1}+0\cdot e_r, \]
and that $\ell=\langle e_1+e_2, \rho\rangle=2r-3$.
Thus, we have
\begin{align*}
\frac{|\rho|^2}{\ell}&=\frac{(r-1)^2+(r-2)^2+\cdots+1^2}{2r-3}=\frac{r(r-1)(2r-1)}{6(2r-3)} > \frac{r^2}{6}=\frac{n}6.
\end{align*}

\item {Type $\mathrm{E}_6$}: For the real split Lie group $G=\mathrm{E}_6^6$, we see from Table \ref{tab:root_data_real_split} that 
\[ \rho=e_2+2e_3+3e_4+4e_5-4e_6-4e_7+4e_8, \]
and that
\[\ell=\langle \frac{1}2(e_1+e_2+e_3+e_4+e_5-e_6-e_7+e_8), \rho\rangle=11.\]
Thus, we have
\begin{align*}
\frac{|\rho|^2}{\ell}&=\frac{1^2+2^2+3^2+4^2+4^2+4^2+4^2}{11} =\frac{78}{11}>7=\frac{n}{6}.
\end{align*}

\item {Type $\mathrm{E}_7$}: For the real split Lie group $G=\mathrm{E}_7^7$, we see from Table \ref{tab:root_data_real_split} that 
        $$\rho=e_2+2e_3+3e_4+4e_5+5e_6-\frac{17}{2}e_7+\frac{17}{2}e_8$$
        and that $\ell=\langle e_8-e_7, \rho\rangle=17$. Thus, we have
        \begin{align*}
            \frac{|\rho|^2}{\ell}&=\frac{1^2+2^2+3^2+4^2+5^2+(17/2)^2+(17/2)^2}{17}=\frac{399}{34}>\frac{35}{3}=\frac{n}{6}.
        \end{align*}
        
\item {Type $\mathrm{E}_8$}: For the real split Lie group $G=\mathrm{E}_8^8$, we see from Table \ref{tab:root_data_real_split} that 
        $$\rho=e_2+2e_3+3e_4+4e_5+5e_6+6e_7+23e_8$$
        and that
        $\ell=\langle e_7+e_8, \rho\rangle=29$. Thus, we have
        \begin{align*}
            \frac{|\rho|^2}{\ell}&=\frac{1^2+2^2+3^2+4^2+5^2+6^2+23^2}{29}=\frac{620}{29}>\frac{64}{3}=\frac{n}{6}.
        \end{align*}
        
\item {Type $\mathrm{F}_4$}: For the real split Lie group $G=\mathrm{F}_4^4$, we see from Table \ref{tab:root_data_real_split} that 
        $$\rho=\frac{11e_1+5e_2+3e_3+e_4}2$$
        and that
        $\ell=\langle e_1+e_2, \rho\rangle=8$. Thus, we have
        \begin{align*}
            \frac{|\rho|^2}{\ell}&=\frac{11^2+5^2+3^2+1^2}{4\times 8}=\frac{39}{8}>\frac{14}{3}=\frac{n}{6}.
        \end{align*}
        
\item {Type $\mathrm{G}_2$}: For the real split Lie group $G=\mathrm{G}_2^2$, we see from Table \ref{tab:root_data_real_split} that
        $$\rho=-e_1-2e_2+3e_3$$
        and that
        $\ell=\langle -e_1-e_2+2e_3, \rho\rangle=9$. Thus, we have
        \begin{align*}
            \frac{|\rho|^2}{\ell}&=\frac{1^2+2^2+3^2}{9}=\frac{14}{9}>\frac{4}{3}=\frac{n}{6}.
        \end{align*}
\end{enumerate}

This exhausts all real split cases. For the remaining cases of each type, the same inequality holds, and we leave the computations in Appendix \ref{app:allcases}.

\appendix

\section{Proof of Theorem \ref{thm:j_X_bound} and Theorem \ref{thm:anosov}: Non-split cases}\label{app:allcases}

The real split cases are already proved in Theorem \ref{thm:j_X_bound} and Theorem \ref{thm:anosov}. We now show that the non-split cases have the same bounds for each type of $G$. The proof is presented in cases. 

\subsection{Type \texorpdfstring{$\mathrm{A}_r$}{Ar}:} 

The corresponding root data (including the set of simple roots together with their multiplicities, the set of all positive roots, the sum of all positive roots counting multiplicities, the half sum of the roots in a maximal strongly orthogonal system, and the dimension of the corresponding symmetric space) is listed in Table \ref{tab:Type_A}

\begin{table}
        \centering
        \begin{tabular}{|c|c|c|}
        \hline
            $G$ & $\SL(r+1, \mathbf{C})$ & $\SU^*(2r+2)$ \\
           \hline
            $(\mathfrak g,\mathfrak a)$ & $\mathrm{A}_r$ & same \\
            \hline
            $\{\alpha_1,\dots, \alpha_r\}$& $\{e_1-e_2,e_2-e_3, \dots, e_r-e_{r+1}\}$ & same \\
            \hline
            Multiplicities & $2,2,\dots, 2$ & $4,4,\dots, 4$\\
            \hline
            $\Delta^+$ & $\{e_i-e_j|i<j\}$ & same \\
            \hline
            $2\rho$ & $(2r)e_1+(2r-4)e_2+\dots+(-2r)e_{r+1}$ & $(4r)e_1+(4r-8)e_2+\dots+(-4r)e_{r+1}$\\
            \hline
            $\Theta$ if $r$ is odd & $\sum_{i=1}^{(r-1)/2}\frac{i}{2}\alpha_i+\sum_{i=(r+1)/2}^{r}\frac{r-i+1}{2}\alpha_i$ & same
            \\
            \hline
            $\Theta$ if $r$ is even & $\sum_{i=1}^{r/2}\frac{i}{2}\alpha_i+\frac{r}{4}\alpha_{r/2+1}+\sum_{i=r/2+2}^{r}\frac{r-i+1}{2}\alpha_i$ & same\\
            \hline
            $n$ & $r(r+2)$ & $r(2r+3)$ \\
            \hline
        \end{tabular}
        \caption{Root data of type $\mathrm{A}_r$}
        \label{tab:Type_A}
\end{table}  

\subsubsection{\texorpdfstring{$G=\SL(r+1,\mathbf{C})$}{G=SL(r+1,C)}:}

{\bf Proof of Theorem \ref{thm:j_X_bound}:} When $r=2m+1$ is odd, we compute that
\[ \Theta=\frac{1}{2}e_1+\dots+\frac{1}{2}e_{m+1}-\frac{1}{2}e_{m+2}-\dots-\frac{1}{2}e_{r+1} \]
and
\[2\rho-\Theta=\frac{8m+3}{2}e_1+\frac{8m-5}{2}e_2+\dots+\frac{3}{2}e_{m+1}-\frac{3}{2}e_{m+2}-\dots-\frac{8m-5}{2}e_r-\frac{8m+3}{2}e_{r+1}.\]
Thus we have
\begin{align*}
\langle \Theta, 2\rho-\Theta\rangle&=2\left(\frac{8m+3}{4}+\frac{8m-5}{4}+\dots+\frac{3}{4}\right)=\frac{(4m+3)(m+1)}{2},
\end{align*}
and $\ell=\langle e_1-e_{r+1}, 2\rho-\Theta\rangle=8m+3$. Therefore, by \eqref{eq:j_X_bound}
\[ j_X(\Gamma)\leq n+1-\frac{\langle \Theta, 2\rho-\Theta\rangle}{\ell} < n+1-\frac{r}{8}.\]
When $r=2m$ is even, we have
\[ \Theta=\frac{1}{2}e_1+\dots+\frac{1}{2}e_{m}+0\cdot e_{m+1}-\frac{1}{2}e_{m+2}-\dots-\frac{1}{2}e_{r+1} \]
and
\[2\rho-\Theta=\frac{8m-1}{2}e_1+\frac{8m-9}{2}e_2+\dots+\frac{7}{2}e_{m}+0\cdot e_{m+1}-\frac{7}{2}e_{m+2}-\dots-\frac{8m-9}{2}e_r-\frac{8m-1}{2}e_{r+1}.\]
Thus we have
\begin{align*}
\langle \Theta, 2\rho-\Theta\rangle&=\left(\frac{8m-1}{4}+\frac{8m-9}{4}+\dots+\frac{7}{4}\right)\times 2=\frac{(4m+3)m}{2},
\end{align*}
and $\ell=\langle e_1-e_{r+1}, 2\rho-\Theta\rangle=8m-1$. Therefore, by \eqref{eq:j_X_bound},
\[j_X(\Gamma)\leq n+1-\frac{\langle \Theta, 2\rho-\Theta\rangle}{\ell} < n+1-\frac{r}{8}.\]
    
{\bf Proof of Theorem \ref{thm:anosov}:} We see from table that
        $$\rho=r\cdot e_1+(r-2)\cdot e_2+\dots+(-r)\cdot e_{r+1},$$
        and that
        $\ell=\langle e_1-e_{r+1}, \rho\rangle=2r$. Thus, when $r=2m+1$ is odd, we have
        \begin{align*}
            \frac{|\rho|^2}{\ell}&=\frac{\left[(2m+1)^2+(2m-1)^2+\cdots+1^2\right]\cdot 2}{2(2m+1)}=\frac{2m^2+5m+3}{3}>\frac{(2m+1)(2m+3)}{6}=\frac{n}6.
        \end{align*}
     when $r=2m$ is even, we have
        \begin{align*}
            \frac{|\rho|^2}{\ell}&=\frac{\left[(2m)^2+(2m-2)^2+\cdots+2^2\right]\cdot 2}{4m}=\frac{2m^2+3m+1}{3}>\frac{2m^2+2m}{3}=\frac{n}6.  
        \end{align*}
        
\subsubsection{\texorpdfstring{$G=\SU^*(2r+2)$}{G=SU*(2r+2)}:}

{\bf Proof of Theorem \ref{thm:j_X_bound}:} When $r=2m+1$ is odd, we compute that
$$\Theta=\frac{1}{2}e_1+\dots+\frac{1}{2}e_{m+1}-\frac{1}{2}e_{m+2}-\dots-\frac{1}{2}e_{r+1}$$
        and
        \[2\rho-\Theta=\frac{16m+7}{2}e_1+\frac{16m-9}{2}e_2+\dots+\frac{7}{2}e_{m+1}-\frac{7}{2}e_{m+2}-\dots-\frac{16m-9}{2}e_r-\frac{16m+7}{2}e_{r+1}.\]
        Thus we have
    \begin{align*}
        \langle \Theta, 2\rho-\Theta\rangle&=\left(\frac{16m+7}{4}+\frac{16m-9}{4}+\dots+\frac{7}{4}\right)\times 2=\frac{(8m+7)(m+1)}{2},
    \end{align*}
    and $\ell=\langle e_1-e_{r+1}, 2\rho-\Theta\rangle=16m+7$. Therefore, by \eqref{eq:j_X_bound},
    \[j_X(\Gamma)\leq n+1-\frac{\langle \Theta, 2\rho-\Theta\rangle}{\ell} < n+1-\frac{r}{8}.\]
    When $r=2m$ is even, we have
    $$\Theta=\frac{1}{2}e_1+\dots+\frac{1}{2}e_{m}+0\cdot e_{m+1}-\frac{1}{2}e_{m+2}-\dots-\frac{1}{2}e_{r+1}$$
        and
        \[2\rho-\Theta=\frac{16m-1}{2}e_1+\frac{16m-17}{2}e_2+\dots+\frac{15}{2}e_{m}+0\cdot e_{m+1}-\frac{15}{2}e_{m+2}-\dots-\frac{16m-17}{2}e_r-\frac{16m-1}{2}e_{r+1}.\]
        Thus we have
    \begin{align*}
        \langle \Theta, 2\rho-\Theta\rangle&=2\left(\frac{16m-1}{4}+\frac{16m-17}{4}+\dots+\frac{15}{4}\right)=\frac{(8m+7)m}{2},
    \end{align*}
    and $\ell=\langle e_1-e_{r+1}, 2\rho-\Theta\rangle=16m-1$. Therefore, by \eqref{eq:j_X_bound},
    \[j_X(\Gamma)\leq n+1-\frac{\langle \Theta, 2\rho-\Theta\rangle}{\ell} < n+1-\frac{r}{8}.\]
    
{\bf Proof of Theorem \ref{thm:anosov}:} We see from table that
        $$\rho=2r\cdot e_1+(2r-4)\cdot e_2+\dots+(-2r)\cdot e_{r+1},$$
        and that
        $\ell=\langle e_1-e_{r+1}, \rho\rangle=4r$. Thus, when $r=2m+1$ is odd, we have
        \begin{align*}
            \frac{|\rho|^2}{\ell}&=\frac{8\left[(2m+1)^2+(2m-1)^2+\cdots+1^2\right]}{4(2m+1)}=\frac{4m^2+10m+6}{3}>\frac{(2m+1)(4m+5)}{6}=\frac{n}6.
        \end{align*}
     when $r=2m$ is even, we have
        \begin{align*}
            \frac{|\rho|^2}{\ell}&=\frac{8\left[(2m)^2+(2m-2)^2+\cdots+2^2\right]}{8m}=\frac{4m^2+6m+2}{3}>\frac{4m^2+3m}{3}=\frac{n}6.  
        \end{align*}

\subsubsection{\texorpdfstring{$G=\mathrm{E}_6^{-26}$}{G=E6(-26)}:}
This Lie group is of type $\mathrm{A}_2$. For Theorem \ref{thm:j_X_bound}, this was already proved in low rank cases $r\leq 3$. For Theorem \ref{thm:anosov}, this is excluded by the condition.

\subsection{Type \texorpdfstring{$\mathrm{B}_r$}{Br}:} The corresponding root data are listed in Table \ref{tab:Type_B}.
\begin{table}
        \centering
        \begin{tabular}{|c|c|c|}
        \hline
            $G$ & $\SO(r,r+k)$ & $\SO(2r+1, \mathbf{C})$ \\
           \hline
            $(\mathfrak g,\mathfrak a)$ & $\mathrm{B}_r$ & same \\
            \hline
            $\{\alpha_1,\dots, \alpha_r\}$& $\{e_1-e_2,e_2-e_3, \dots, e_{r-1}-e_r, e_r\}$ & same \\
            \hline
            Multiplicities & $1,1,\dots,1, k$ & $2,2,\dots, 2,2$\\
            \hline
            $\Delta^+$ & $\{e_i\pm e_j|i<j\}\cup \{e_i\}$ & same \\
            \hline
            $2\rho$ & $(2r-2+k)e_1+(2r-4+k)e_2+\dots+(k)e_{r}$ & $(4r-2)e_1+(4r-6)e_2+\dots+2e_{r}$\\
            \hline
            $\Theta$ if $r$ is odd & $\sum_{i=1}^{(r-1)/2}i\alpha_i+\sum_{i=(r+1)/2}^{r}\frac{r}{2}\alpha_i$  & same
            \\
            \hline
            $\Theta$ if $r$ is even &  $\sum_{i=1}^{r/2}i\alpha_i+\sum_{i=r/2+1}^{r}\frac{r}{2}\alpha_i$  & same\\
            \hline
            $n$ & $r(r+k)$ & $r(2r+1)$ \\
            \hline
        \end{tabular}
        \caption{Root data of type $\mathrm{B}_r$}
        \label{tab:Type_B}
\end{table}  

\subsubsection{\texorpdfstring{$G=\SO(r,r+k)$}{}:}

{\bf Proof of Theorem \ref{thm:j_X_bound}:} 

Recall that we already showed the low rank cases $r=2,3$, so we may assume $r\geq 4$. When $r=2m+1$ is odd, we compute that $$\Theta=e_1+\dots+e_{m}+\frac{1}{2}e_{m+1}$$
        and
        \begin{align*}
            2\rho-\Theta=&(4m-1+k)e_1+(4m-3+k)e_2+\dots+(2m+1+k)e_{m}+(2m-\frac{1}{2}+k)e_{m+1}\\
            &+(2m-2+k)e_{m+2}+\dots+(k)e_{2m+1}
        \end{align*}
        
        Thus we have
    \begin{align*}
        \langle \Theta, 2\rho-\Theta\rangle&=(4m-1+k)+(4m-3+k)+\dots+(2m+1+k)+\frac{4m-1+2k}{4}\\
        &=3m^2+km+m+\frac{k}{2}-\frac{1}{4}
    \end{align*}
    and $\ell=\langle e_1+e_2, 2\rho-\Theta\rangle=8m-4+2k$. It follows that
    \[\frac{\langle \Theta, 2\rho-\Theta\rangle}{\ell}=\frac{3r}{16}+\frac{2mk+8m+k+4}{16(4m+k-2)}>\frac{3r}{16}\]    
    Therefore, by \eqref{eq:j_X_bound},
    \[j_X(\Gamma)\leq n+1-\frac{\langle \Theta, 2\rho-\Theta\rangle}{\ell} < n+1-\frac{3r}{16}.\]

When $r=2m$ is even, we compute that $$\Theta=e_1+\dots+e_{m}$$
        and
        \begin{align*}
            2\rho-\Theta=&(4m-3+k)e_1+(4m-5+k)e_2+\dots+(2m-1+k)e_{m}+(2m-2+k)e_{m+1}\\
            &+(2m-4+k)e_{m+2}+\dots+(k)e_{2m}
        \end{align*}
        Thus we have
    \begin{align*}
        \langle \Theta, 2\rho-\Theta\rangle&=(2r-3+k)+(2r-5+k)+\dots+(r-1+k)=\frac{r(3r-4+2k)}{4},
    \end{align*}
    and $\ell=\langle e_1+e_2, 2\rho-\Theta\rangle=4r-8+2k$. It follows that
    \[\frac{\langle \Theta, 2\rho-\Theta\rangle}{\ell}=\frac{3r}{16}+\frac{rk+4r}{16(2r+k-4)}>\frac{3r}{16}.\]   
    Therefore, by \eqref{eq:j_X_bound},
    \[j_X(\Gamma)\leq n+1-\frac{\langle \Theta, 2\rho-\Theta\rangle}{\ell} < n+1-\frac{3r}{16}.\]
    
\begin{remark}
    In the stable case when $k\rightarrow \infty$, the bound can be improved to $j_X(\Gamma)\leq n+1-\frac{r}{4}$.
\end{remark}

{\bf Proof of Theorem \ref{thm:anosov}:} We see from table that
        $$\rho=(r-1+\frac{k}2)\cdot e_1+(r-2+\frac{k}2)\cdot e_2+\dots+\frac{k}2\cdot e_{r},$$
        and that
        $\ell=\langle e_1+e_{2}, \rho\rangle=2r-3+k$. Thus, we have
        \begin{align*}
            \frac{|\rho|^2}{\ell}&=\frac{(r-1+\frac{k}2)^2+(r-2+\frac{k}2)^2+\cdots+(\frac{k}2)^2}{2r-3+k}\\
            &=\frac{\sum_{i=1}^r (i^2+(k-2)i+(\frac{k}{2}-1)^2)}{2r-3+k}\\
            &=\frac{r}{6}\cdot\frac{2r^2+3rk-3r+\frac{3}{2}k^2-3k+1}{2r-3+k} >\frac{r(r+k)}{6}=\frac{n}6.
        \end{align*}

\subsubsection{\texorpdfstring{$G=\SO(2r+1,\mathbf{C})$}{}:}

{\bf Proof of Theorem \ref{thm:j_X_bound}:} We already showed the low rank cases $r=2,3$, so we may assume $r\geq 4$.  

When $r=2m+1$ is odd, we compute that $$\Theta=e_1+\dots+e_{m}+\frac{1}{2}e_{m+1}$$
        and
\begin{align*}
     2\rho-\Theta=&(8m+1)e_1+(8m-3)e_2+\dots+(4m+5)e_{m}+(4m+\frac{3}{2})e_{m+1}\\
     &+(4m-2)e_{m+2}+\dots+2e_{2m+1}
\end{align*}
        Thus we have
    \begin{align*}
        \langle \Theta, 2\rho-\Theta\rangle&=(8m+1)+(8m-3)+\dots+(4m+5)+(2m+\frac{3}{4})=\frac{6r^2-2r-1}{4},
    \end{align*}
    and $\ell=\langle e_1+e_2, 2\rho-\Theta\rangle=8r-10$. It follows that
    \[\frac{\langle \Theta, 2\rho-\Theta\rangle}{\ell}=\frac{\frac{6r^2-2r-1}{4}}{8r-10}>\frac{3r}{16}.\]    
    Therefore, by \eqref{eq:j_X_bound},
    \[j_X(\Gamma)\leq n+1-\frac{\langle \Theta, 2\rho-\Theta\rangle}{\ell} < n+1-\frac{3r}{16}.\]
When $r=2m$ is even, we have
$$\Theta=e_1+\dots+e_{m}$$
and
\begin{align*}
     2\rho-\Theta=&(8m-3)e_1+(8m-7)e_2+\dots+(4m+1)e_{m}+\\
     &+(4m-2)e_{m+1}+\dots+2e_{2m}
\end{align*}
Thus we have
    \begin{align*}
        \langle \Theta, 2\rho-\Theta\rangle&=(8m-3)+(8m-7)+\dots+(4m+1)=\frac{3r^2-r}{2},
    \end{align*}
    and $\ell=\langle e_1+e_2, 2\rho-\Theta\rangle=8r-10$. It follows that
    \[\frac{\langle \Theta, 2\rho-\Theta\rangle}{\ell}=\frac{\frac{3r^2-r}{2}}{8r-10}>\frac{3r}{16}.\]

    Therefore, by \eqref{eq:j_X_bound},
    \[j_X(\Gamma)\leq n+1-\frac{\langle \Theta, 2\rho-\Theta\rangle}{\ell} < n+1-\frac{3r}{16}.\]
    
{\bf Proof of Theorem \ref{thm:anosov}:} We see from table that
        $$\rho=(2r-1)\cdot e_1+(2r-3)\cdot e_2+\dots+1\cdot e_{r},$$
        and that
        $\ell=\langle e_1+e_{2}, \rho\rangle=4r-4$. Thus, we have
        \begin{align*}
            \frac{|\rho|^2}{\ell}&=\frac{(2r-1)^2+(2r-3)^2+\cdots+1^2}{4r-4}=\frac{r(4r^2-1)}{3(4r-4)}>\frac{r(2r+1)}{6}=\frac{n}6.
        \end{align*}
    
\subsection{Type \texorpdfstring{$\mathrm{C}_r$}{}:}  

The corresponding root data are listed in Table \ref{tab:Type_C}.
\begin{table}[ht]
        \centering
        \begin{tabular}{|c|c|c|}
        \hline
            $G$ & $\SU(r,r)$ & $\Sp(r,\mathbf{C})$ \\
           \hline
            $(\mathfrak g,\mathfrak a)$ & $\mathrm{C}_r$ & same \\
            \hline
            $\{\alpha_1,\dots, \alpha_r\}$& $\{e_1-e_2,e_2-e_3, \dots, e_{r-1}-e_r, 2e_r\}$ & same \\
            \hline
            Multiplicities & $2,2,\dots,2, 1$ & $2,2,\dots, 2,2$\\
            \hline
            $\Delta^+$ & $\{e_i\pm e_j|i<j\}\cup \{2e_i\}$ & same \\
            \hline
            $2\rho$ & $(4r-2)e_1+(4r-6)e_2+\dots+2e_{r}$ & $(4r)e_1+(4r-4)e_2+\dots+4e_{r}$\\
            \hline
            $\Theta$ & $\sum_{i=1}^{r-1}i\alpha_i+\frac{r}{2}\alpha_r$  & same
            \\
            \hline
            $n$ & $2r^2$ & $r(2r+1)$\\
            \hline
        \end{tabular}
        \begin{tabular}{|c|c|c|}
        \hline
            $G$ & $\SO^*(4r)$ & $\Sp(r,r)$ \\
           \hline
            $(\mathfrak g,\mathfrak a)$ & same & same \\
            \hline
            $\{\alpha_1,\dots, \alpha_r\}$& same & same \\
            \hline
            Multiplicities & $4,4,\dots,4, 1$ & $4,4,\dots, 4,3$\\
            \hline
            $\Delta^+$ & same & same \\
            \hline
            $2\rho$ & $(8r-6)e_1+(8r-14)e_2+\dots+2e_{r}$ & $(8r-2)e_1+(8r-10)e_2+\dots+6e_{r}$\\
            \hline
            $\Theta$ & same  & same
            \\
             \hline
            $n$ & $2r(2r-1)$ & $4r^2$\\
            \hline
        \end{tabular}
        \begin{tabular}{|c|c|}
        \hline
            $G$ & $\mathrm{E}_7^{-25}$  \\
           \hline
            $(\mathfrak g,\mathfrak a)$ & $\mathrm{C}_3$ \\
            \hline
            $\{\alpha_1,\alpha_2, \alpha_3\}$& same \\
            \hline
            Multiplicities & $8,8,1$ \\
            \hline
            $\Delta^+$ & same  \\
            \hline
            $2\rho$ & $34e_1+18e_2+2e_3$\\
            \hline
            $\Theta$ & same  \\
             \hline
            $n$ & $54$\\
            \hline
        \end{tabular}
        \caption{Root data of type $\mathrm{C}_r$}
        \label{tab:Type_C}
\end{table}  

\subsubsection{\texorpdfstring{$G=\SU(r,r)$}{}:}

{\bf Proof of Theorem \ref{thm:j_X_bound}:} We compute that $$\Theta=e_1+\dots+e_{r}$$
        and
\begin{align*}
     2\rho-\Theta=&(4r-3)e_1+(4r-7)e_2+\dots+e_{r}.
\end{align*}
        Thus we have
    \begin{align*}
        \langle \Theta, 2\rho-\Theta\rangle&=(4r-3)+(4r-7)+\dots+1=(2r-1)r,
    \end{align*}
    and $\ell=\langle 2e_1, 2\rho-\Theta\rangle=8r-6$. Therefore, by \eqref{eq:j_X_bound},
    \[j_X(\Gamma)\leq n+1-\frac{\langle \Theta, 2\rho-\Theta\rangle}{\ell} < n+1-\frac{r}{4}.\]  
{\bf Proof of Theorem \ref{thm:anosov}:} We see from table that
        $$\rho=(2r-1)\cdot e_1+(2r-3)\cdot e_2+\dots+1\cdot e_{r},$$
        and that
        $\ell=\langle 2e_1, \rho\rangle=4r-2$. Thus, we have
        \begin{align*}
            \frac{|\rho|^2}{\ell}&=\frac{(2r-1)^2+(2r-3)^2+\cdots+1^2}{4r-2}=\frac{r(4r^2-1)}{3(4r-2)}>\frac{r^2}{3}=\frac{n}6.
        \end{align*}

\subsubsection{\texorpdfstring{$G=\Sp(r,\mathbf{C})$}{}:}

{\bf Proof of Theorem \ref{thm:j_X_bound}:} We compute that $$\Theta=e_1+\dots+e_{r}$$
        and
\begin{align*}
     2\rho-\Theta=&(4r-1)e_1+(4r-5)e_2+\dots+3e_{r}.
\end{align*}
        Thus we have
    \begin{align*}
        \langle \Theta, 2\rho-\Theta\rangle&=(4r-1)+(4r-5)+\dots+3=(2r+1)r,
    \end{align*}
    and $\ell=\langle 2e_1, 2\rho-\Theta\rangle=8r-2$. Therefore, by \eqref{eq:j_X_bound},
    \[j_X(\Gamma)\leq n+1-\frac{\langle \Theta, 2\rho-\Theta\rangle}{\ell} < n+1-\frac{r}{4}.\]
{\bf Proof of Theorem \ref{thm:anosov}:} We see from table that
        $$\rho=(2r)\cdot e_1+(2r-2)\cdot e_2+\dots+2\cdot e_{r},$$
        and that
        $\ell=\langle 2e_1, \rho\rangle=4r$. Thus, we have
        \begin{align*}
            \frac{|\rho|^2}{\ell}&=\frac{(2r)^2+(2r-2)^2+\cdots+2^2}{4r}=\frac{(r+1)(2r+1)}{6}>\frac{r(2r+1)}{6}=\frac{n}6.
        \end{align*}
\subsubsection{\texorpdfstring{$G=\SO^*(4r)$}{}:}

{\bf Proof of Theorem \ref{thm:j_X_bound}:} We compute that $$\Theta=e_1+\dots+e_{r}$$
        and
\begin{align*}
     2\rho-\Theta=&(8r-7)e_1+(8r-15)e_2+\dots+e_{r}.
\end{align*}
        Thus we have
    \begin{align*}
        \langle \Theta, 2\rho-\Theta\rangle&=(8r-7)+(8r-15)+\dots+1=(4r-3)r,
    \end{align*}
    and $\ell=\langle 2e_1, 2\rho-\Theta\rangle=16r-14$. Therefore, by \eqref{eq:j_X_bound},
    \[j_X(\Gamma)\leq n+1-\frac{\langle \Theta, 2\rho-\Theta\rangle}{\ell} < n+1-\frac{r}{4}.\]
{\bf Proof of Theorem \ref{thm:anosov}:} We see from table that
        $$\rho=(4r-3)\cdot e_1+(4r-7)\cdot e_2+\dots+1\cdot e_{r},$$
        and that
        $\ell=\langle 2e_1, \rho\rangle=8r-6$. Thus, we have
        \begin{align*}
            \frac{|\rho|^2}{\ell}&=\frac{(4r-3)^2+(4r-7)^2+\cdots+1^2}{8r-6}=\frac{r(16r^2-12r-1)}{6(4r-3)}>\frac{r(2r-1)}{3}=\frac{n}6.
        \end{align*}

\subsubsection{\texorpdfstring{$G=\Sp(r,r)$}{}:}

{\bf Proof of Theorem \ref{thm:j_X_bound}:} We compute that $$\Theta=e_1+\dots+e_{r}$$
        and
\begin{align*}
     2\rho-\Theta=&(8r-3)e_1+(8r-11)e_2+\dots+5e_{r}.
\end{align*}
        Thus we have
    \begin{align*}
        \langle \Theta, 2\rho-\Theta\rangle&=(8r-3)+(8r-11)+\dots+5=(4r+1)r,
    \end{align*}
    and $\ell=\langle 2e_1, 2\rho-\Theta\rangle=16r-6$. Therefore, by \eqref{eq:j_X_bound},
    \[j_X(\Gamma)\leq n+1-\frac{\langle \Theta, 2\rho-\Theta\rangle}{\ell} < n+1-\frac{r}{4}.\]
{\bf Proof of Theorem \ref{thm:anosov}:} We see from table that
        $$\rho=(4r-1)\cdot e_1+(4r-5)\cdot e_2+\dots+3\cdot e_{r},$$
        and that
        $\ell=\langle 2e_1, \rho\rangle=8r-2$. Thus, we have
        \begin{align*}
            \frac{|\rho|^2}{\ell}&=\frac{(4r-1)^2+(4r-5)^2+\cdots+3^2}{8r-2}=\frac{r(16r^2+12r-1)}{6(4r-1)}>\frac{2r^2}{3}=\frac{n}6.
        \end{align*}

\subsubsection{\texorpdfstring{$G=\mathrm{E}_7^{-25}$}{}:} 

{\bf Proof of Theorem \ref{thm:j_X_bound}:} The Lie group is of type $\mathrm{C}_3$, which was already proved in low rank cases $r\leq 3$.

{\bf Proof of Theorem \ref{thm:anosov}:} We see from table that
        $$\rho=17e_1+9e_2+e_3,$$
        and that
        $\ell=\langle 2e_1, \rho\rangle=34$. Thus, we have
        \begin{align*}
            \frac{|\rho|^2}{\ell}&=\frac{17^2+9^2+1^2}{34}=\frac{371}{34}>9=\frac{n}6.
        \end{align*}
\subsection{Type \texorpdfstring{$\mathrm{D}_r$}{}:} 

The only non-real split example is $G=\SO(2r,\mathbf{C})$. The corresponding root data are listed in Table \ref{tab:Type_D}. 
\begin{table}
        \centering
        \begin{tabular}{|c|c|}
        \hline
            $G$ & $\SO(2r,\mathbf{C})$ \\
           \hline
            $(\mathfrak g,\mathfrak a)$ & $\mathrm{D}_r$  \\
            \hline
            $\{\alpha_1,\dots, \alpha_r\}$& $\{e_1-e_2,e_2-e_3, \dots, e_{r-1}-e_r, e_{r-1}+e_r\}$ \\
            \hline
            Multiplicities & $2,2,\dots,2, 2$ \\
            \hline
            $\Delta^+$ & $\{e_i\pm e_j|i<j\}$  \\
            \hline
            $2\rho$ & $(4r-4)e_1+(4r-8)e_2+\dots+4e_{r-1}$ \\
            \hline
             $\Theta$ if $r$ is odd  &$\sum_{i=1}^{(r-1)/2}i\alpha_i+\sum_{i=(r+1)/2}^{r-2}\frac{r-1}{2}\alpha_i+\frac{r-1}{4}(\alpha_{r-1}+\alpha_r)$  \\
            \hline
            $\Theta$ if $r$ is even 
            & $\sum_{i=1}^{r/2}i\alpha_i+\sum_{i=r/2+1}^{r-2}\frac{r}{2}\alpha_i+\frac{r}{4}(\alpha_{r-1}+\alpha_r)$ \\
            \hline
            $n$ & $r(2r-1)$ \\
            \hline
        \end{tabular}
        \caption{Root data of type $\mathrm{D}_r$}
        \label{tab:Type_D}
\end{table}  

\subsubsection{\texorpdfstring{$G=\SO(2r,\mathbf{C})$}{}:}

{\bf Proof of Theorem \ref{thm:j_X_bound}:} When $r=2m+1$ is odd, we compute that $$\Theta=e_1+\dots+e_{m}$$
        and
\begin{align*}
     2\rho-\Theta=&(8m-1)e_1+(8m-5)e_2+\dots+(4m+3)e_{m}\\
     &+(4m)e_{m+1}+(4m-4)e_{m+2}+\dots+4e_{2m}\\
\end{align*}
        Thus we have
    \begin{align*}
        \langle \Theta, 2\rho-\Theta\rangle&=(8m-1)+(8m-5)+\dots+(4m+3)=\frac{3r^2-5r+2}{2},
    \end{align*}
    and $\ell=\langle e_1+e_2, 2\rho-\Theta\rangle=8r-14$. Therefore, by \eqref{eq:j_X_bound},
    \[j_X(\Gamma)\leq n+1-\frac{\langle \Theta, 2\rho-\Theta\rangle}{\ell} < n+1-\frac{3r}{16}.\]
When $r=2m$ is even, we have
$$\Theta=e_1+\dots+e_{m}$$
and
\begin{align*}
     2\rho-\Theta=&(8m-5)e_1+(8m-9)e_2+\dots+(4m-1)e_{m}+\\
     &+(4m-4)e_{m+1}+(4m-8)e_{m+2}+\dots+4e_{2m-1}
\end{align*}
Thus we have
    \begin{align*}
        \langle \Theta, 2\rho-\Theta\rangle&=(8m-5)+(8m-9)+\dots+(4m-1)=\frac{3r(r-1)}{2},
    \end{align*}
    and $\ell=\langle e_1+e_2, 2\rho-\Theta\rangle=8r-14$. Therefore, by \eqref{eq:j_X_bound},
    \[j_X(\Gamma)\leq n+1-\frac{\langle \Theta, 2\rho-\Theta\rangle}{\ell} < n+1-\frac{3r}{16}.\]

{\bf Proof of Theorem \ref{thm:anosov}:} We see from table that
        $$\rho=(2r-2)\cdot e_1+(2r-4)\cdot e_2+\dots+2\cdot e_{r-1},$$
        and that
        $\ell=\langle e_1+e_2, \rho\rangle=4r-6$. Thus, we have
        \begin{align*}
            \frac{|\rho|^2}{\ell}&=\frac{(2r-2)^2+(2r-4)^2+\cdots+2^2}{4r-6}=\frac{r(r-1)(2r-1)}{3(2r-3)}>\frac{r(2r-1)}{6}=\frac{n}6.
        \end{align*}
        
\subsection{Type \texorpdfstring{$(\mathrm{BC})_r$}{}:}

The corresponding root data are listed in Table \ref{tab:Type_BC}. We note from \cite{Oh02} that $\Theta$ is defined using the set of all non-multipliable roots, so it is same as in the case of $\mathrm{C}_r$ type. More precisely, if we choose the simple roots as $\alpha_1,\dots,(\alpha_r,2\alpha_r)$, then $\Theta=\sum_{i=1}^{r-1} i\alpha_i+r\alpha_r$. (Note that the coefficient in $\alpha_r$ differs by half to keep the notation in Table \ref{tab:Type_C})

\begin{table}
        \centering
        \begin{tabular}{|c|c|}
        \hline
             $G$ & $\SU(r,r+k)$  \\
           \hline
            $(\mathfrak g,\mathfrak a)$ & $(\mathrm{BC})_r$ \\
            \hline
            $\{\alpha_1,\dots, (\alpha_r,2\alpha_r)\}$& $\{e_1-e_2,e_2-e_3, \dots, e_{r-1}-e_r, (e_r, 2e_r)\}$\\
            \hline
            Multiplicities & $2,2,\dots,2, (2k,1)$ \\
            \hline
            $\Delta^+$ & $\{e_i\pm e_j|i<j\}\cup \{e_i, 2e_i\}$  \\
            \hline
            $2\rho$ & $(4r+2k-2)e_1+(4r+2k-6)e_2+\dots+(2k+2)e_{r}$ \\
            \hline
            $\Theta$ & $\sum_{i=1}^{r}i\alpha_i$
            \\
            \hline
            $n$ & $2r(r+k)$ \\
            \hline
        \end{tabular}
        \begin{tabular}{|c|c|}
        \hline
            $G$ & $\Sp(r,r+k)$  \\
           \hline
            Multiplicities & $4,4,\dots,4, (4k,3)$ \\
            \hline
            $2\rho$ & $(8r+4k-2)e_1+(8r+4k-10)e_2+\dots+(4k+6)e_{r}$ \\
            \hline
            $n$ & $4r(r+k)$ \\
            \hline
            Any other data & same 
            \\
            \hline
        \end{tabular}
        \begin{tabular}{|c|c|c|}
        \hline
            $G$ & $\SO^*(4r+2)$ & $\mathrm{E}_6^{-14}$\\
           \hline
            Multiplicities & $4,4,\dots,4, (4,1)$ & $6,(8,1)$\\
            \hline
            $2\rho$ & $(8r-2)e_1+(8r-10)e_2+\dots+6e_{r}$ & $22e_1+10e_2$\\
            \hline
            $n$ & $2r(2r+1)$ & 32\\
            \hline
            Any other data & same & same
            \\
            \hline
        \end{tabular}
        \caption{Root data of type $(\mathrm{BC})_r$}
        \label{tab:Type_BC}
\end{table}

\subsubsection{\texorpdfstring{$G=\SU(r,r+k)$}{}:}

{\bf Proof of Theorem \ref{thm:j_X_bound}:} We compute that $$\Theta=e_1+\dots+e_{r}$$
        and
\begin{align*}
     2\rho-\Theta=&(4r+2k-3)e_1+(4r+2k-7)e_2+\dots+(2k+1)e_{r}.
\end{align*}
        Thus we have
    \begin{align*}
        \langle \Theta, 2\rho-\Theta\rangle&=(4r+2k-3)+(4r+2k-7)+\dots+(2k+1)=(2r+2k-1)r,
    \end{align*}
    and $\ell=\langle 2e_1, 2\rho-\Theta\rangle=8r+4k-6$. It follows that
    \[\frac{\langle \Theta, 2\rho-\Theta\rangle}{\ell}=\frac{r}{4}+\frac{2kr+r}{4(4r+2k-3)}>\frac{r}{4}\]    
    Therefore, by \eqref{eq:j_X_bound},
    \[j_X(\Gamma)\leq n+1-\frac{\langle \Theta, 2\rho-\Theta\rangle}{\ell} < n+1-\frac{r}{4}.\]
\begin{remark}
    In the stable case when $k\rightarrow \infty$, the bound can be improved to $j_X(\Gamma)\leq n+1-\frac{r}{2}$.
\end{remark}

{\bf Proof of Theorem \ref{thm:anosov}:} We see from table that
        $$\rho=(2r+k-1)\cdot e_1+(2r+k-3)\cdot e_2+\dots+(k+1)\cdot e_{r},$$
        and that
        $\ell=\langle 2e_1, \rho\rangle=4r+2k-2$. Thus, we have
        \begin{align*}
            \frac{|\rho|^2}{\ell}&=\frac{(2r+k-1)^2+(2r+k-3)^2+\cdots+(k+1)^2}{4r+2k-2}\\
            &=\frac{\sum_{i=1}^r(4i^2+4(k-1)i+(k-1)^2)}{4r+2k-2}\\
            &=\frac{r}{3}\left(\frac{4r^2+6kr+3k^2-1}{4r+2k-2}\right)>\frac{r(r+k)}{3}=\frac{n}6.
        \end{align*}

\subsubsection{\texorpdfstring{$G=\Sp(r,r+k)$}{}:}

{\bf Proof of Theorem \ref{thm:j_X_bound}:} We compute that $$\Theta=e_1+\dots+e_{r}$$
        and
\begin{align*}
     2\rho-\Theta=&(8r+4k-3)e_1+(8r+4k-11)e_2+\dots+(4k+5)e_{r}.
\end{align*}
        Thus we have
    \begin{align*}
        \langle \Theta, 2\rho-\Theta\rangle&=(8r+4k-3)+(8r+4k-11)+\dots+(4k+5)=(4r+4k+1)r,
    \end{align*}
    and $\ell=\langle 2e_1, 2\rho-\Theta\rangle=16r+8k-6$. It follows that
    \[\frac{\langle \Theta, 2\rho-\Theta\rangle}{\ell}=\frac{r}{4}+\frac{4kr+5r}{4(8r+4k-3)}>\frac{r}{4}\]    
    Therefore, by \eqref{eq:j_X_bound},
    \[j_X(\Gamma)\leq n+1-\frac{\langle \Theta, 2\rho-\Theta\rangle}{\ell} < n+1-\frac{r}{4}.\]
    
\begin{remark}
    In the stable case when $k\rightarrow \infty$, the bound can be improved to $j_X(\Gamma)\leq n+1-\frac{r}{2}$.
\end{remark}

{\bf Proof of Theorem \ref{thm:anosov}:} We see from table that
        $$\rho=(4r+2k-1)\cdot e_1+(4r+2k-5)\cdot e_2+\dots+(2k+3)\cdot e_{r},$$
        and that
        $\ell=\langle 2e_1, \rho\rangle=2(4r+2k-1)$. Thus, we have
        \begin{align*}
            \frac{|\rho|^2}{\ell}&=\frac{(4r+2k-1)^2+(4r+2k-5)^2+\cdots+(2k+3)^2}{2(4r+2k-1)}\\
            &=\frac{\sum_{i=1}^r(16i^2+8(2k-1)i+(2k-1)^2)}{2(4r+2k-1)}\\
            &=\frac{r}{3}\left(\frac{16r^2+24kr+12r+12k^2+12k-1}{2(4r+2k-1)}\right)>\frac{2r(r+k)}{3}=\frac{n}6.
        \end{align*}

\subsubsection{\texorpdfstring{$G=\SO^*(4r+2)$}{}:}

{\bf Proof of Theorem \ref{thm:j_X_bound}:} We compute that $$\Theta=e_1+\dots+e_{r}$$
        and
\begin{align*}
     2\rho-\Theta=&(8r-3)e_1+(8r-11)e_2+\dots+5e_{r}.
\end{align*}
        Thus we have
    \begin{align*}
        \langle \Theta, 2\rho-\Theta\rangle&=(8r-3)+(8r-11)+\dots+5=(4r+1)r,
    \end{align*}
    and $\ell=\langle 2e_1, 2\rho-\Theta\rangle=16r-6$. Therefore, by \eqref{eq:j_X_bound},
    \[j_X(\Gamma)\leq n+1-\frac{\langle \Theta, 2\rho-\Theta\rangle}{\ell} < n+1-\frac{r}{4}.\]
{\bf Proof of Theorem \ref{thm:anosov}:} We see from table that
        $$\rho=(4r-1)\cdot e_1+(4r-5)\cdot e_2+\dots+3\cdot e_{r},$$
        and that
        $\ell=\langle 2e_1, \rho\rangle=8r-2$. Thus, we have
        \begin{align*}
            \frac{|\rho|^2}{\ell}&=\frac{(4r-1)^2+(4r-5)^2+\cdots+3^2}{8r-2}=\frac{r(16r^2+12r-1)}{6(4r-1)}>\frac{r(2r+1)}{3}=\frac{n}6.
        \end{align*}

\subsubsection{\texorpdfstring{$G=\mathrm{E}_6^{-14}$}{}:}
{\bf Proof of Theorem \ref{thm:j_X_bound}:} The Lie group is of type $(\mathrm{BC})_2$, which was already proved in low rank cases $r\leq 3$.

{\bf Proof of Theorem \ref{thm:anosov}:} We see from table that
        $$\rho=11e_1+5e_2,$$
        and that
        $\ell=\langle 2e_1, \rho\rangle=22$. Thus, we have
        \begin{align*}
            \frac{|\rho|^2}{\ell}&=\frac{11^2+5^2}{22}=\frac{73}{11}>\frac{16}{3}=\frac{n}6.
        \end{align*}
\subsection{Type \texorpdfstring{$\mathrm{E}_6$}{}:}The only non-real split example is $G=\mathrm{E}_6(\mathbf{C})$. The corresponding root data are listed in Table \ref{tab:E}.

\begin{table}
        \centering
        \begin{tabular}{|c|c|}
        \hline
            $G$ & $\mathrm{E}_6(\mathbf{C})$ \\
           \hline
            $(\mathfrak g,\mathfrak a)$ & $E_6$  \\
            \hline
           $\{\alpha_1,\dots, \alpha_6\}$& $\{\frac{1}{2}(e_8-e_7-e_6-e_5-e_4-e_3-e_2+e_1), e_2+e_1, e_2-e_1, e_3-e_2,e_4-e_3,e_5-e_4\}$\\
            \hline
            Multiplicities & $2,2,2,2,2,2$ \\
            \hline
            $\Delta^+$ & $\{e_i\pm e_j|5\geq i>j\}\cup \{\frac{1}{2}(e_8-e_7-e_6+\sum_{i=1}^5(-1)^{n(i)}e_i)|\sum_i n(i) \textrm{ is even}\}$ \\
            \hline
            $2\rho$ & $4e_2+8e_3+12e_4+16e_5-16e_6-16e_7+16e_8$\\
            \hline
            $n$ & $78$\\
            \hline
        \end{tabular}
         \begin{tabular}{|c|c|}
        \hline
            $G$ & $\mathrm{E}_7(\mathbf{C})$ \\
           \hline
            $(\mathfrak g,\mathfrak a)$ & $\mathrm{E}_7$  \\
            \hline
           $\{\alpha_1,\dots, \alpha_7\}$& $\{\frac{1}{2}(e_8-e_7-e_6-e_5-e_4-e_3-e_2+e_1), e_2+e_1\}\cup\{e_{i+1}-e_i|1\leq i\leq 5\}\}$\\
            \hline
            Multiplicities & $2,2,2,2,2,2,2$ \\
            \hline
           $\Delta^+$ & $\{e_i\pm e_j|6\geq i>j\}\cup\{e_8-e_7\}\cup \{\frac{1}{2}(e_8-e_7+\sum_{i=1}^6(-1)^{n(i)}e_i)|\sum_i n(i) \textrm{ is odd}\}$ \\
            \hline
            $2\rho$ & $4e_2+8e_3+12e_4+16e_5+20e_6-34e_7+34e_8$\\
            \hline
            $n$ & $133$\\
            \hline
        \end{tabular}
         \begin{tabular}{|c|c|}
        \hline
            $G$ & $\mathrm{E}_8(\mathbf{C})$ \\
           \hline
            $(\mathfrak g,\mathfrak a)$ & $\mathrm{E}_8$  \\
            \hline
           $\{\alpha_1,\dots, \alpha_8\}$& $\{\frac{1}{2}(e_8-e_7-e_6-e_5-e_4-e_3-e_2+e_1), e_2+e_1\}\cup\{e_{i+1}-e_i|1\leq i\leq 6\}$\\
            \hline
            Multiplicities & $2,2,2,2,2,2,2,2$ \\
            \hline
           $\Delta^+$ & $\{e_i\pm e_j|i>j\}\cup \{\frac{1}{2}(e_8+\sum_{i=1}^7(-1)^{n(i)}e_i)|\sum_i n(i) \textrm{ is even}\}$ \\
            \hline
            $2\rho$ & $4e_2+8e_3+12e_4+16e_5+20e_6+24e_7+92e_8$\\
            \hline
            $n$ & $248$\\
            \hline
        \end{tabular}
        \caption{Root data of type $\mathrm{E}_6, \mathrm{E}_7, \mathrm{E}_8$}
        \label{tab:E}
\end{table}  

\subsubsection{\texorpdfstring{$G=\mathrm{E}_6(\mathbf{C})$}{}:}

{\bf Proof of Theorem \ref{thm:anosov}:} We see from table that
        $$\rho=2e_2+4e_3+6e_4+8e_5-8e_6-8e_7+8e_8,$$
        and that
        \[\ell=\langle \frac{1}2(e_1+e_2+e_3+e_4+e_5-e_6-e_7+e_8), \rho\rangle=22.\]
        Thus, we have
        \begin{align*}
            \frac{|\rho|^2}{\ell}&=\frac{2^2+4^2+6^2+8^2\times 4}{22}=\frac{156}{11}>13=\frac{n}6.
        \end{align*}

\subsection{Type \texorpdfstring{$\mathrm{E}_7$}{}:}

The only non-real split example is $G=\mathrm{E}_7(\mathbf{C})$. The corresponding root data are listed in Table \ref{tab:E}. 
\subsubsection{\texorpdfstring{$G=\mathrm{E}_7(\mathbf{C})$}{}:}
{\bf Proof of Theorem \ref{thm:anosov}:} We see from table that
        $$\rho=2e_2+4e_3+6e_4+8e_5+10e_6-17e_7+17e_8,$$
        and that
        $\ell=\langle e_8-e_7, \rho\rangle=34$. Thus, we have
        \begin{align*}
            \frac{|\rho|^2}{\ell}&=\frac{2^2+4^2+6^2+8^2+10^2+17^2+17^2}{34}=\frac{399}{17}>\frac{133}{6}=\frac{n}6.
        \end{align*}
        
\subsection{Type \texorpdfstring{$\mathrm{E}_8$}{}:}

The only non-real split example is $G=\mathrm{E}_8(\mathbf{C})$. The corresponding root data are listed in Table \ref{tab:E}.
\subsubsection{\texorpdfstring{$G=\mathrm{E}_8(\mathbf{C})$}{}:}
{\bf Proof of Theorem \ref{thm:anosov}:} We see from table that
        $$\rho=2e_2+4e_3+6e_4+8e_5+10e_6+12e_7+46e_8,$$
        and that
        $\ell=\langle e_8+e_7, \rho\rangle=58$. Thus, we have
        \begin{align*}
            \frac{|\rho|^2}{\ell}&=\frac{2^2+4^2+6^2+8^2+10^2+12^2+46^2}{58}=\frac{1240}{29}>\frac{124}{3}=\frac{n}6.
        \end{align*}

\subsection{Type \texorpdfstring{$\mathrm{F}_4$}{}:} 

The corresponding root data are listed in Table \ref{tab:F}.
\begin{table}
        \centering
        \begin{tabular}{|c|c|c|}
        \hline
            $G$ & $\mathrm{E}_6^2$ & $\mathrm{E}_7^{-5}$ \\
           \hline
            $(\mathfrak g,\mathfrak a)$ & $\mathrm{F}_4$ & same \\
            \hline
           $\{\alpha_1,\dots, \alpha_4\}$& $\{\frac{1}{2}(e_1-e_2-e_3-e_4), e_4, e_3-e_4, e_2-e_3\}$ & same\\
            \hline
            Multiplicities & $2,2,1,1$ & $4,4,1,1$\\
            \hline
             $\Delta^+$ & $\{e_i\pm e_j|i<j\}\cup \{e_i\}\cup \{\frac{1}{2}(e_1\pm e_2\pm e_3\pm e_4)\}$ & same\\
            \hline
            $2\rho$ & $16e_1+6e_2+4e_3+2e_4$ & $26e_1+8e_2+6e_3+4e_4$ \\
            \hline
            $n$ & $40$ & $64$\\
            \hline
        \end{tabular}
        \begin{tabular}{|c|c|c|}
        \hline
            $G$ & $\mathrm{E}_8^{-24}$ & $\mathrm{F}_4(\mathbf{C})$\\
           \hline
            Multiplicities & $8,8,1,1$ & $2,2,2,2$\\
            \hline
            $2\rho$ & $46e_1+12e_2+10e_3+8e_4$& $22e_1+10e_2+6e_3+2e_4$\\
            \hline
            $n$ & $112$ & $52$\\
            \hline
            Any other data & same & same
            \\
            \hline
        \end{tabular}
        \begin{tabular}{|c|c|}
           \hline
            $(\mathfrak g, \mathfrak a)$ & $\mathrm{G}_2$  \\
            \hline
            $\{\alpha_1,\alpha_2\}$& $\{e_1-e_2, -2e_1+e_2+e_3\}$\\
            \hline
            Multiplicities & $2,2$ \\
            \hline
            $\Delta^+$ & $\{e_1-e_2, e_3-e_1, e_3-e_2,-2e_1+e_2+e_3, e_1-2e_2+e_3, -e_1-e_2+2e_3\}$ \\
            \hline
            $2\rho$ & $-4e_1-8e_2+12e_3$\\
            \hline
             $n$ & $14$\\
            \hline
        \end{tabular}
        \caption{Root data of type $\mathrm{F}_4$ and $\mathrm{G}_2$}
        \label{tab:F}
\end{table}

\subsubsection{\texorpdfstring{$G=\mathrm{E}_6^2$}{}:}

{\bf Proof of Theorem \ref{thm:anosov}:} We see from table that
        $$\rho=8e_1+3e_2+2e_3+e_4,$$
        and that
        $\ell=\langle e_1+e_2, \rho\rangle=11$. Thus, we have
        \begin{align*}
            \frac{|\rho|^2}{\ell}&=\frac{8^2+3^2+2^2+1^2}{11}=\frac{78}{11}>\frac{20}{3}=\frac{n}6.
        \end{align*}
        
\subsubsection{\texorpdfstring{$G=\mathrm{E}_7^{-5}$}{}:}

{\bf Proof of Theorem \ref{thm:anosov}:} We see from table that
        $$\rho=13e_1+4e_2+3e_3+2e_4,$$
        and that
        $\ell=\langle e_1+e_2, \rho\rangle=17$. Thus, we have
        \begin{align*}
            \frac{|\rho|^2}{\ell}&=\frac{13^2+4^2+3^2+2^2}{17}=\frac{198}{17}>\frac{32}{3}=\frac{n}6.
        \end{align*}
        
\subsubsection{\texorpdfstring{$G=\mathrm{E}_8^{-24}$}{}:}

{\bf Proof of Theorem \ref{thm:anosov}:} We see from table that
        $$\rho=23e_1+6e_2+5e_3+4e_4,$$
        and that
        $\ell=\langle e_1+e_2, \rho\rangle=29$. Thus, we have
        \begin{align*}
            \frac{|\rho|^2}{\ell}&=\frac{23^2+6^2+5^2+4^2}{29}=\frac{606}{29}>\frac{56}{3}=\frac{n}6.
        \end{align*}
        
\subsubsection{\texorpdfstring{$G=\mathrm{F}_4(\mathbf{C})$}{}:}

{\bf Proof of Theorem \ref{thm:anosov}:} We see from table that
        $$\rho=11e_1+5e_2+3e_3+e_4,$$
        and that
        $\ell=\langle e_1+e_2, \rho\rangle=16$. Thus, we have
        \begin{align*}
            \frac{|\rho|^2}{\ell}&=\frac{11^2+5^2+3^2+1^2}{16}=\frac{39}{4}>\frac{26}{3}=\frac{n}6.
        \end{align*}

\subsection{Type \texorpdfstring{$\mathrm{G}_2$}{}:}

The only non-real split example is $G=\mathrm{G}_2(\mathbf{C})$. The corresponding root data are listed in Table \ref{tab:F}.

\subsubsection{\texorpdfstring{$G=\mathrm{G}_2(\mathbf{C})$}{}:}

{\bf Proof of Theorem \ref{thm:anosov}:} We see from table that
        $$\rho=-2e_1-4e_2+6e_3,$$
        and that
        $\ell=\langle -e_1-e_2+2e_3, \rho\rangle=18$. Thus, we have
        \begin{align*}
            \frac{|\rho|^2}{\ell}&=\frac{2^2+4^2+6^2}{18}=\frac{28}{9}>\frac{7}{3}=\frac{n}6.
        \end{align*}

This exhausts all non-real split cases, hence the proofs of Theorem \ref{thm:j_X_bound} and Theorem \ref{thm:anosov} are now complete.



\quad

CC: Indiana University, Bloomington, IN, USA. E-mail: \verb|connell@indiana.edu|

DBM: Purdue University, West Lafayette, IN, USA. E-mail: \verb|dmcreyno@purdue.edu|

SW: ShanghaiTech University, Pudong, Shanghai, China. E-mail: \verb|wangshi@shanghaitech.edu.cn|


\end{document}